\documentclass[1p,preprint]{elsarticle}
\journal{Linear Algebra and its Applications}
\usepackage{amsmath,amssymb,todonotes}
\usepackage{exscale,theorem}
\usepackage{verbatim,graphicx,xcolor}
 \usepackage{tikz}
\usetikzlibrary{arrows}
 \usepackage{hyperref}
\usepackage{float}
\usepackage{subfigure}
\usepackage{url}

\usepackage[color]{changebar}
\cbcolor{blue}
\usepackage{lipsum}

\newcommand{\ignore}[1]{}

\usepackage[normalem]{ulem}

\DeclareMathAlphabet{\mathpzc}{OT1}{pzc}{m}{it}
\newtheorem{theorem}{Theorem}[section]

\newtheorem{corollary}[theorem]{Corollary}
\newtheorem{proposition}[theorem]{Proposition}
\newtheorem{definition}[theorem]{Definition}

{\theorembodyfont{\rm}}
{\theorembodyfont{\rm}\newtheorem{remark}[subsection]{Remark}}

\begin{document}
\begin{frontmatter}

\title{Eigenvalue avoidance of structured matrices depending smoothly on a real parameter\tnoteref{t1}}

\author[nii]{Yuji Nakatsukasa\fnref{label1}}\ead{nakatsukasa@maths.ox.ac.uk}

\author[ess]{Vanni Noferini\corref{cor1}\fnref{label3}}
\ead{vanni.noferini@aalto.fi}

\address[nii]{
Mathematical Institute, University of Oxford, Oxford, OX2 6GG}

\address[ess]{Department of Mathematics and Systems Analysis,
Aalto University, PO Box 11100, 00076 Aalto, Finland.
Supported by an Academy of Finland grant (Suomen Akatemian p\"{a}\"{a}tos 331240).}
\cortext[cor1]{Corresponding author.}

\setlength{\marginparsep}{0.25cm}
\setlength{\marginparwidth}{1.5cm}
\def\mnote#1{{%
              \marginpar{\color{red}\flushright\raggedright\footnotesize\textsf{#1}}}}

\newenvironment{vn}{\begin{quote}\color{violet} \small\sf VN $\diamondsuit$~}{\end{quote}}
\newenvironment{yn}{\begin{quote}\color{dgreen} \small\sf YN $\diamondsuit$~}{\end{quote}}

\begin{abstract}
We explore the concept of eigenvalue avoidance, which is well understood for real symmetric and Hermitian matrices, for other classes of structured matrices. We adopt a differential geometric perspective and study the generic behaviour of the eigenvalues of regular and injective curves $t \in ]a,b[ \mapsto A(t) \in \mathcal{N} $ where $\mathcal{N}$ is a smooth real Riemannian submanifold of either $\mathbb{R}^{n \times n}$ or $\mathbb{C}^{n \times n}$. We focus on the case where $\mathcal{N}$ corresponds to some class of (real or complex) structured matrices including skew-symmetric, skew-Hermitian, orthogonal, unitary, banded symmetric, banded Hermitian, banded skew-symmetric, and banded skew-Hermitian. We argue that for some structures eigenvalue avoidance always happens, whereas for other structures this may depend on the parity of the size, on the numerical value of the multiple eigenvalue, and possibly on the value of determinant. As a further application of our tools we also study singular value avoidance for unstructured matrices in $\mathbb{R}^{m \times n}$ or $\mathbb{C}^{m \times n}$. 
\end{abstract}


\end{frontmatter}

\makeatletter
\def\ps@pprintTitle{%
     \let\@oddhead\@empty
     \let\@evenhead\@empty
     \def\@oddfoot{}
     \let\@evenfoot\@oddfoot}
\makeatother

\newproof{proof}{Proof}
\let\scshape\bfseries  

\def\qedsymbol{\vbox{\hrule\hbox{%
                     \vrule height1.3ex\hskip0.8ex\vrule}\hrule}}
\def\endproof{\qquad\qedsymbol\medskip\par}
\def\myendproof{\qquad\qedsymbol}
\def\noqed{\def\qedsymbol{}}

\def\leapproxeq{\mathrel{\raisebox{-.75ex}{$\mathop{\sim}\limits^{\textstyle <}$}}}
\def\geapproxeq{\mathrel{\raisebox{-.75ex}{$\mathop{\sim}\limits^{\textstyle >}$}}}

\newcommand {\mpar}[1]{\marginpar{\fussy\tiny #1}}
\def\R{\mathbb{R}}
\def\C{\mathbb{C}}
\def\G{\mathbb{G}}
\def\Gt{\widetilde{\mathbb{G}}}
\def\Ah{\widehat{A}}
\def\L{\mathbb{L}}
\def\K{\mathbb{F}}
\def\J{\mathbb{J}}
\def\SS{\mathbb{S}}
\def\R{\mathbb{R}}
\def\C{\mathbb{C}}
\def\cS{\mathcal{S}}
\def\cH{\mathcal{H}}
\def\cW{\mathcal{W}}
\def\cM{\mathcal{M}}
\def\cN{\mathcal{N}}
\def\N{\mathbb{N}}
\def\U{\mathcal{U}}
\def\V{\mathcal{V}}
\def\re{\mathop{\mathrm{Re}}}
\def\im{\mathop{\mathrm{Im}}}
\def\F{Fr\'{e}chet }
\def\Im{\mathrm{Im}}
\def\nbyn{n \times n}
\def\mbyn{m \times n}
\def\l{\lambda}
\def\normt#1{\|#1\|_2}
\def\normf#1{\|#1\|_F}
\def\norm#1{\|#1\|}
\def\normi#1{\|#1\|_1}
\def\normo#1{\|#1\|_{\infty}}
\def\Chat{\widehat{C}}
\def\e{eigenvalue}
\def\diag{\mathop{\mathrm{diag}}}     
\def\trace{\mathop{\mathrm{trace}}}
\def\sign{\mathop{\mathrm{sign}}}
\def\sig{\mathop{\mathrm{sig}}}
\def\err{\mathop{\mathrm{err}}}
\def\At{\widetilde{A}}
\def\normt#1{\|#1\|_2}
\def\struc{\textup{struc}}
\definecolor{dgreen}{rgb}{0.13, 0.55, 0.13}

\def\condabs{\textup{cond}_\textup{abs}}
\def\vec{\textup{vec}}
\def\unvec{\textup{unvec}}
\def\Adj{\textup{adj}}
\def\addots{\mathinner{\mkern1mu\raise1pt\hbox{.}\mkern2mu\raise4pt\hbox{.}
     \mkern2mu\raise7pt\vbox{\kern7pt\hbox{.}}\mkern1mu}}

\newcommand{\pth}[1]{\textit{p}th #1}
\newcommand\scalemath[2]{\scalebox{#1}{\mbox{\ensuremath{\displaystyle #2}}}}
\def\sgn{\mathrm{\mathop{sign}}}
\def\rank{\mathrm{\mathop{rank}}}
\def\range{\mathrm{\mathop{range}}}
\def\adj{\mathrm{\mathop{adj}}}
\def\Null{\mathrm{\mathop{null}}}
\def\spann{\mathrm{\mathop{span}}}
\def\scond{\mathrm{\mathop{stcond}}}
\def\cond{\mathrm{\mathop{cond}}}
\def\ubcond{\mathrm{\mathop{ub}}_-\mathrm{\mathop{cond}}}

\def\lbcond{\mathrm{\mathop{lb}}_-\mathrm{\mathop{cond}}}

\def\mystrut#1{\rule{0cm}{#1}}

\makeatletter
\def\mymatrix#1{\null\,\vcenter{\normalbaselines\m@th
    \ialign{\hfil$##$\hfil&&\quad\hfil$##$\hfil\crcr
      \mathstrut\crcr\noalign{\kern-\baselineskip}
      #1\crcr\mathstrut\crcr\noalign{\kern-\baselineskip}}}\,}
\makeatother
\def\mybmatrix#1{\left[ \mymatrix{#1} \right]}
\def\twobytwo#1#2#3#4{\bigl[{\hfil#1\atop\hfil#3}{\hfil#2\atop\hfil#4}\bigr]}
\def\dtwobytwo#1#2#3#4{\bigl|{\hfil#1\atop\hfil#3}{\hfil#2\atop\hfil#4}\bigr|}
\def\twobyone#1#2{\bigl[{\hfil#1\atop\hfil#2}\bigr]}

\def\scalar#1{\langle #1\rangle}
\def\sM{{{}_{\mbox {\tiny M}}}}

\def\sMthm{{{}_{\mbox {\tiny {\em M}}}}}
\def\sIJthm{{{}_{\mbox {\tiny {\em IJ}}}}}
\def\sJ{{}_{\mbox {\tiny J}}}
\def\sS{{}_{\mbox {\tiny S}}}
\def\BSM{B^{}_{\SS_{\sMthm}}}
\def\BGM{B^{}_{\G_{\sMthm}}}
\def\BLM{B^{}_{\L_{\sMthm}}}
\def\BGIJ{B^{}_{\G_{\sIJthm}}}
\def\BGt{B^{}_{\Gt}}

\makeatletter
\def\revddots{\mathinner{\mkern1mu\raise\p@
    \vbox{\kern7\p@\hbox{.}}\mkern2mu
    \raise4\p@\hbox{.}\mkern2mu\raise7\p@\hbox{.}\mkern1mu}}
\makeatother


\mathcode`@="8000 
 {\catcode`\@=\active\gdef@{\mkern1mu}}
\newcounter{mylineno}
\makeatletter
\let\oldtabcr\@tabcr
\def\nonumberbreak{\oldtabcr\hspace{3.5pt}}
\def\mynewline{\refstepcounter{mylineno}%
                \llap{\footnotesize\arabic{mylineno}\hspace{5pt}}%
               }
\def\lineref#1{\footnotesize\ref{#1}}
\gdef\@tabcr{\@stopline \@ifstar{\penalty%
            \@M \@xtabcr}\@xtabcr\mynewline}
\def\myvspace#1{\oldtabcr[#1]\mynewline}
\newenvironment{code}{%
                         \mathcode`\:="603A  
                         \def\colon{\mathchar"303A}
                         \setcounter{mylineno}{0}
                         \par
                         \upshape
                         \begin{list} 
                            {} {\leftmargin = 1cm}
                         \item[]
                         \begin{tabbing}

                            \hspace*{.3in} \= \hspace*{.3in} \=
                            \hspace*{.3in} \= \hspace*{.3in} \= \kill
                            \mynewline
                        }{\end{tabbing}\end{list}}
\makeatother

\setlength{\marginparsep}{0.25cm}
\setlength{\marginparwidth}{1.5cm}
\def\mnote#1{{%
              \marginpar{\color{red}\flushright\raggedright\footnotesize\textsf{#1}}}}

\newenvironment{djh}{\begin{quote}\color{blue} \small\sf DJH $\diamondsuit$~}{\end{quote}}
\newenvironment{vn}{\begin{quote}\color{violet} \small\sf VN $\diamondsuit$~}{\end{quote}}
\newenvironment{pg}{\begin{quote}\color{dgreen} \small\sf PG $\diamondsuit$~}{\end{quote}}

\section{Introduction}
Eigenvalue avoidance for self-adjoint 
 operators depending smoothly on a real parameter $t$ is a phenomenon originally
 observed by physicists while developing quantum mechanics. 
It was first explained by von Neumann and Wigner~\cite{von1993verhalten} arguing on the codimension of the set of real symmetric (or Hermitian) matrices with multiple eigenvalues. 
The theory was generalized to infinite-dimensional self-adjoint operators 
by Uhlenbeck~\cite{uhlenbeck1976generic} and Teytel~\cite{teytel1999rare}. 
This phenomenon is related to 
various applications; see~\cite{betcke2004computations} 
 for more information and background.

 In this work we study eigenvalue avoidance for the case of structured finite dimensional linear operators, i.e., matrices. We consider regular \cite{smale1958} injective curves of matrices that exhibit some kind of structure, different than being symmetric or Hermitian.  We take a geometric approach and examine the real codimension of the submanifold of structured matrices with multiple eigenvalues, when embedded in the real manifold of matrices of the same structure (with general eigenvalues). All the structures that we consider in our study of eigenvalue avoidance are special cases of normal matrices. For a normal matrix, having  a multiple eigenvalue is equivalent to being a derogatory matrix \cite[Definitions 1.4.4]{hornjohn}, i.e., a matrix with an eigenvalue whose geometric multiplicity is higher than $1$. If the codimension is 2 or more, and if $t \mapsto A(t)$ is a 
regular and injective curve on the manifold of the structured matrices under consideration, then  generically there is no value of $t$ such that $A(t)$ has multiple eigenvalues; a fact that is a consequence of classical results in differential geometry, as we discuss in Section~\ref{sec:geom}. 
We investigate a variety of matrix structures, namely orthogonal/unitary, skew-symmetric, and banded (skew-)symmetric. We also examine the singular values of real and complex matrices. In Table~\ref{tab:summary} we summarize  our findings along with previous results for $n\times n$ matrices.

\begin{table}[htbp]
  \centering
  \caption{Ambient (real) dimension and 
codimension of derogatory matrices for various matrix structures. 
}
  \label{tab:summary}
  \begin{tabular}{c|c|c|c}

 & ambient dimension&  codimension &  reference \\\hline\hline
symmetric & $\frac{n(n+1)}{2}$ & 2 & \cite{Laxlaa,von1993verhalten} \\
Hermitian & $n^2$ & 3  & \cite{keller2008,Laxlaa}\\\hline
skew-symmetric, $n$ even & $\frac{n(n+1)}{2}$ & 1  & Thm.~\ref{thm:skewsym}\\
skew-symmetric, $n$ odd & $\frac{n(n+1)}{2}$ & 3  & Thm.~\ref{thm:skewsym} \\
skew-Hermitian & $n^2$ & 3  & Thm.~\ref{thm:skewHerm}\\\hline
orthogonal, $n$ even, $\det(A)=1$ & $\frac{n(n-1)}{2}$ & 1  & Thm.~\ref{thm:orth} \\
orthogonal, $n$ even, $\det(A)=-1$ & $\frac{n(n-1)}{2}$ & 3  & Thm.~\ref{thm:orth} \\
orthogonal, $n$ odd & $\frac{n(n-1)}{2}$ & 1  & Thm.~\ref{thm:orth} \\
unitary & $n^2$ & 3  & \cite{keller2008} \\\hline
banded & (same as dense) & (same as dense)  & \cite{DPP18} \\\hline
singular values $\mathbb{R}^{m\times n}$& $mn$ & 2  & \cite{keller2008} \\
singular values $\mathbb{C}^{m\times n}$& $2mn$ & 3  & Thm.~\ref{thm:svd} \\
  \end{tabular}
\end{table}

Table~\ref{tab:summary} also provides older references for those structures for which we were able to find some treatment in the literature; in the case of structures for which no reference is reported, we may be the first to provide an analysis. Hence, our paper is in part a review that proposes a new point of view on known facts. On the other hand, we also provide results that are (we believe) novel and interesting, in particular the codimensions of skew-symmetric matrices depend crucially on the parity of $n$, and the location of the eigenvalue collision (zero or not). Similarly, the codimension of orthogonal matrices with multiple eigenvalues depends on an intricate combination of the parity and the determinant.

Before starting our journey through derogatory 
structured matrices and their codimensions, it is appropriate to give credit to a number of mathematicians who led the way. As previously mentioned, Von Neumann and Wigner opened this line of research in \cite{von1993verhalten} by studying the case of real symmetric matrices. More recently, Keller~\cite{keller2008} examined the codimension of matrices with prescribed eigenvalue multiplicities, treating Hermitian, normal, and (complex) unitary matrices. Keller also considered the singular values of real matrices. 
The structures that are considered in this paper, but not in \cite{keller2008}, are real orthogonal matrices, skew-symmetric  and skew-Hermitian matrices, their banded versions, and singular values of complex matrices. 
We show that interesting phenomena arise in each of these cases. 

We also treat the cases already studied in \cite{keller2008}. The line of argument used by Keller for counting the codimension is similar to ours, but not the same; moreover, our focus is on the implications in terms of eigenvalue avoidance. Even more recently, for banded matrices, Dieci, Papini and Pugliese~\cite{DPP18} mentioned that the codimensions are the same as the full counterparts due to previous results in~\cite{dieci1999smooth}. 
However, not many details are given, and (at least to us) the implication is not immediate; for this reason, as well as for the sake of self-containedness, we prefer to include our own treatment.

 \section{A geometric perspective on eigenvalue avoidance}\label{sec:geom}

A simple experiment can be set up to observe eigenvalue avoidance. If $A(t)$ is an $n\times n$ real symmetric matrix whose entries are continuous (resp. analytic) functions of $t \in ]a,b[$ then there are $n$ continuous (resp. analytic)  functions $\lambda_1(t),\lambda_2(t),\dots,\lambda_n(t)$ that, for each $t \in ]a,b[$, are the eigenvalues of $A(t)$ (see~\cite[Sec. II.5]{kato} and~\cite[Sec. 5]{MNTX16} for more details and \citep{rellich37,rellichbook} as historical references). For example, one could simply let $A(t)=A_0 + A_1 t$ for some $A_0=A_0^T, A_1=A_1^T \in \R^{n \times n}$; for instance, one could generate $A_0$ and $A_1$ from some random matrix ensemble. What typically happens is that eigenvalue paths may come very close to each other, but they do not cross; Figure~\ref{fig1} below gives an example. Of course, it is possible to construct $A_0$ and $A_1$ in such a way that a crossing does happen, but this appears to be rare event in some sense. Von Neumann and Wigner \cite{von1993verhalten} were the first to explain this lack of crossing with a geometric argument. To explain what the connection is, we start by recalling the definition of codimension of an embededd submanifold.

\begin{definition}
Let $\cM$ be a smooth submanifold of the smooth manifold $\cN$. Then the codimension of $\cM$ when embededd in $\cN$ is
$$\dim \cN - \dim \cM.$$
\end{definition}

\begin{figure}[htbp]
  \begin{minipage}[t]{0.33\hsize}
\includegraphics[width=.98\textwidth]{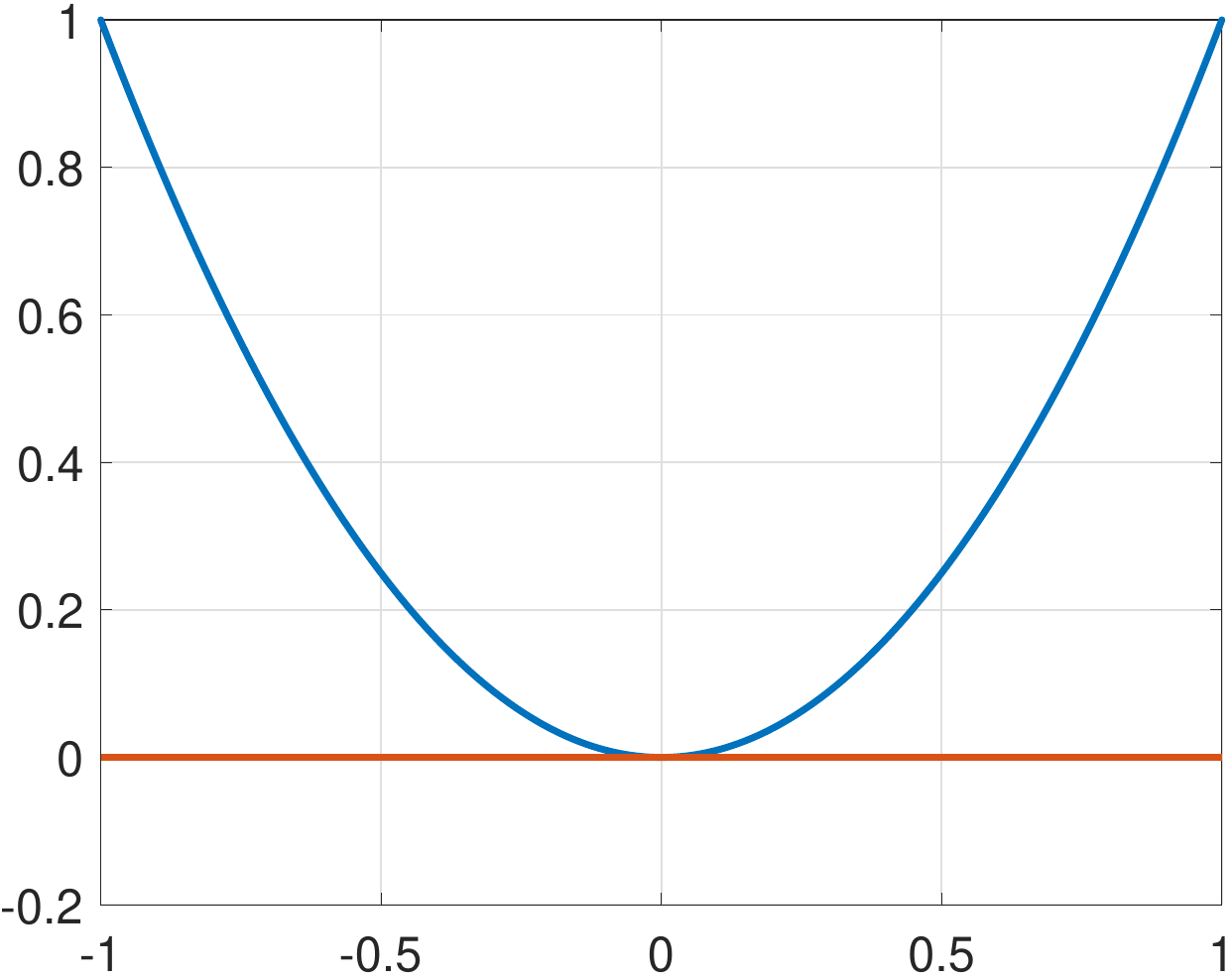} 
  \end{minipage}   
  \begin{minipage}[t]{0.325\hsize}
\includegraphics[width=1\textwidth]{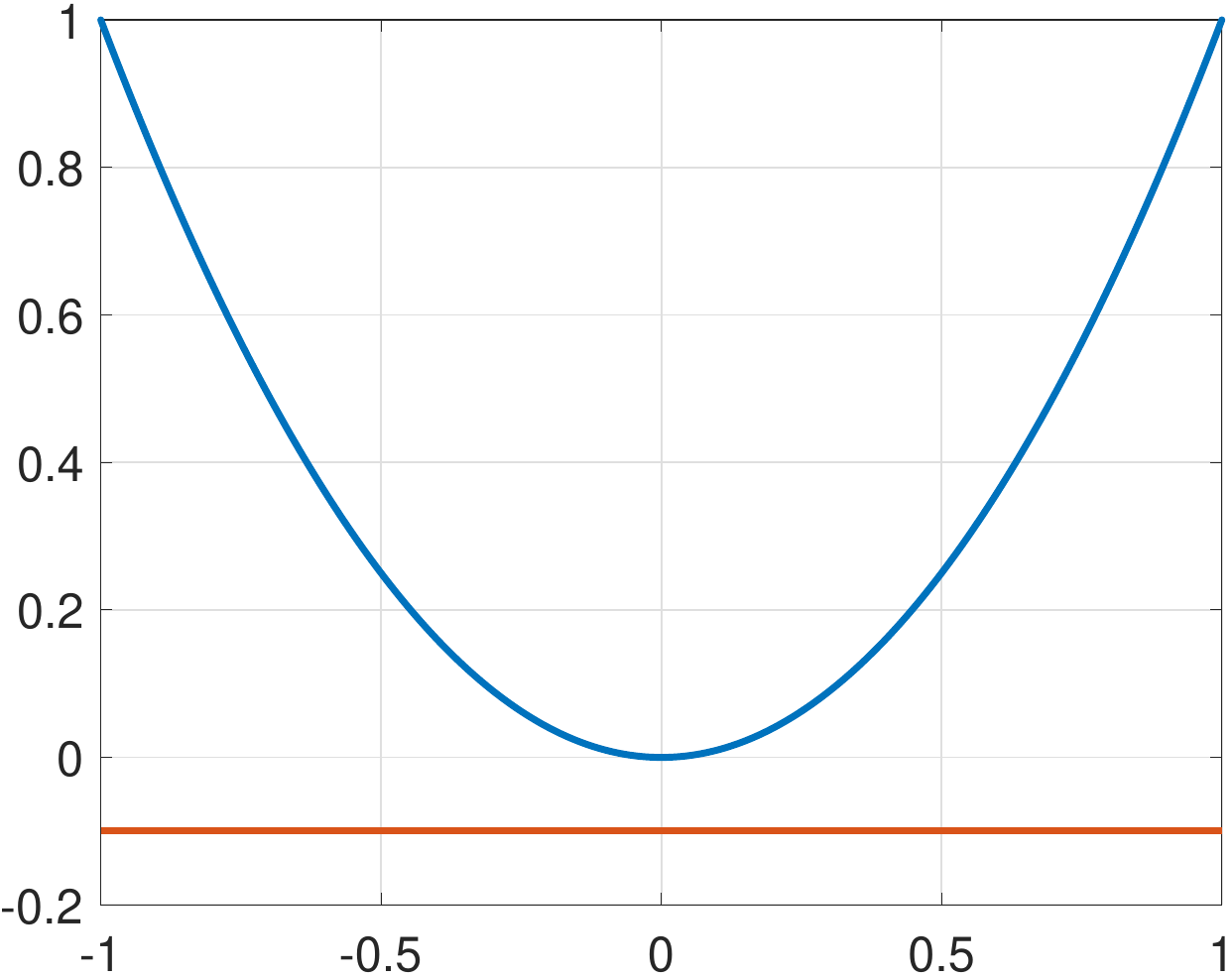}
  \end{minipage}
  \begin{minipage}[t]{0.33\hsize}
\includegraphics[width=.98\textwidth]{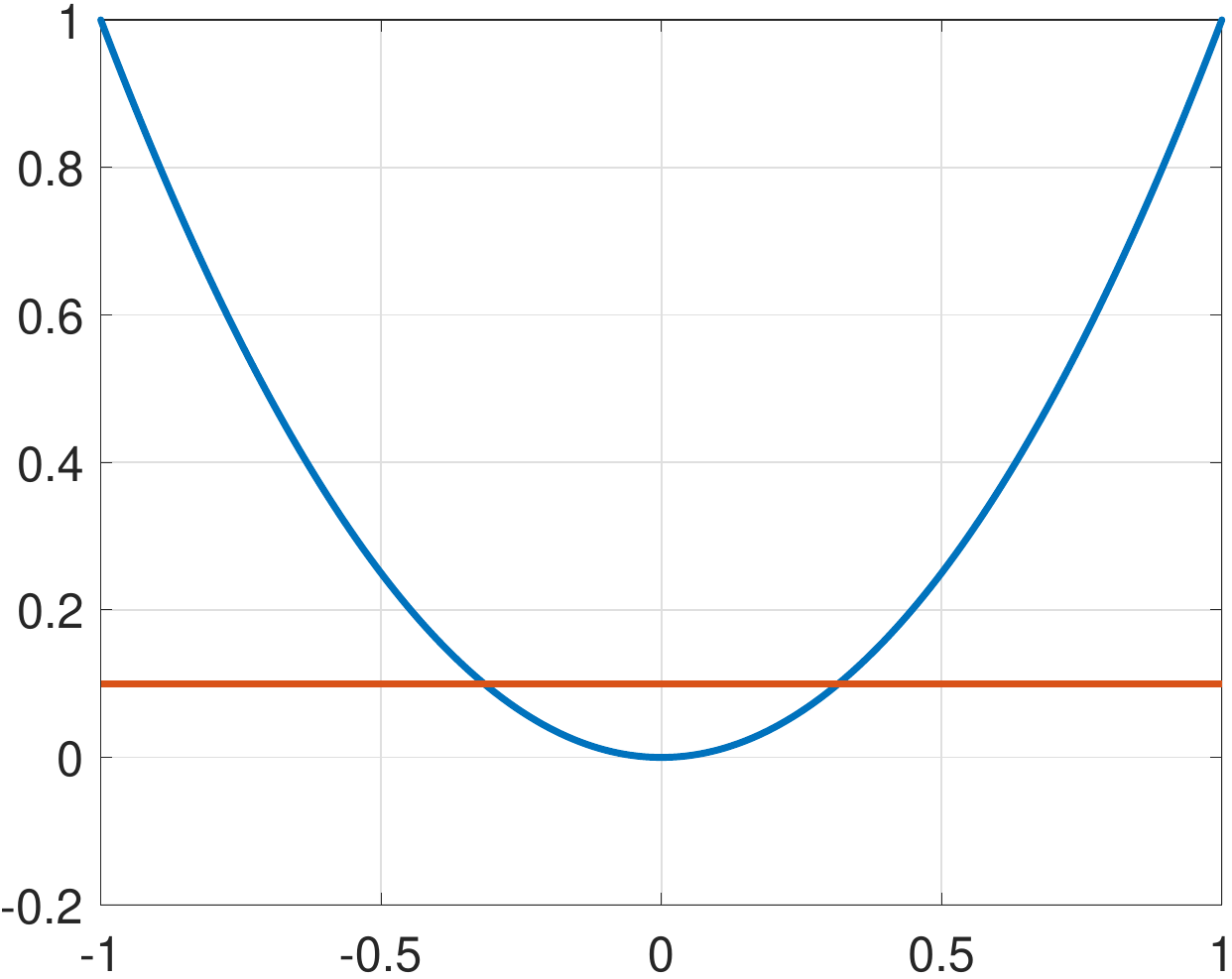}
  \end{minipage}
\caption{ 
Two submanifolds of $\R^2$, a parabula (blue) and a line (red), that happen not to inteserct transversally (left). By slightly translating the red line, the two manifolds either intersect transversally (right) or do not intersect at all and hence are transversal by definition (center).}
\label{fig00}\end{figure}

Codimension is a particularly useful concept in intersection theory, because it is generically additive: recall that two submanifolds $\cM_0$ and $\cM_1$ are said to be \emph{transversal}, or to intersect transversally, if the vector space sum of their tangent spaces at any point in the intersection generates the tangent space of the whole $\cN$ at that point \cite[Sec. 6]{lee}. We note  the following important subtlety: since the previous definition is of the form ``for all points in the intersection $\dots$", it implies that any two submanifolds that do \emph{not} intersect are transversal. If $\cM_0$ and $\cM_1$ are transversal, then the codimension of their intersection is the sum of their codimensions: 
$$\dim \cN - \dim (\cM_0 \cap \cM_1) = 2 \dim \cN - \dim \cM_0 - \dim \cM_1.$$
Not all manifolds intersect transversally, but transversality is generic. 
This can be made precise, for example, as follows: fix two real smooth submanifolds of $V \in \{\R^{n \times n}, \C^{n \times n}\}$, say, $\mathcal{X}$ and $\mathcal{Y}$, and let $x \in V$. Then, for almost all (in the Euclidean topology) $x$, $x + \mathcal{X}$ and $\mathcal{Y}$ intersect transversally. Or, to put it in a simpler but perhaps more vivid way: If you are so unlucky to have found two manifolds that \emph{do not} intersect transversally, imagine trying to fix the situation by slightly translating one of them; for almost all such small tilts, the two manifolds will then intersect transversally. See Figure \ref{fig00}.

The statement that transversality of intersections is generic holds beyond the special case of translations; we have chosen to focus on it simply because it is, in our opinion, illustrative. A more complete treatment of this subject is beyond the scope of this paper, but more details can be found, for example, in \cite[Sec. 6]{lee}. For us, it suffices to say that these results on genericity stem from a celebrated result in differential geometry known as Thom's Transversality Theorem \cite{arnold2012geometrical,Thom2,Thom}. In particular, we have
\begin{corollary}\label{cor:thom}
Let $\cM_0$, $\cM_1$ be submanifolds of $\cN$.
If $\dim \cM_0 + \dim \cM_1 < \dim \cN$, then generically $\cM_0$ and $\cM_1$ do not intersect.
\end{corollary}
\begin{proof}
The only way for $\cM_0$ and $\cM_1$ to be (vacuously) transversal is not to intersect at all.
\end{proof}
Suppose that $\cN$ is the manifold of real symmetric matrices, and $\cM_0$ is the submanifold of derogatory real symmetric matrices. Von Neumann and Wigner showed that the codimension of $\cM_0$, when embedded in $\cN$, is $2$. Now, if $A(t)$ is a regular and injective curve, and hence a smooth submanifold of dimension $1$, by Corollary \ref{cor:thom} we conclude that it generically will not intersect $\cM_0$, i.e., it will not have multiple eigenvalues.

\section{Counting codimensions of derogatory structured matrices}

In the next sections we expose the dimensional argument of Von Neumann and Wigner that, via Corollary \ref{cor:thom}, explains eigenvalue avoidance of real symmetric and Hermitian matrices continuous in $t$. Then, we generalize it to predict eigenvalue avoidance (or, in some cases, lack thereof) for other structured matrices. More specifically, for real matrices $X \in \R^{n \times n}$ we consider the following structures:
\begin{itemize}
\item symmetric matrices, $X=X^T$;
\item skew-symmetric matrices, $X=-X^T$;
\item orthogonal matrices, $X X^T= I_n$;
\item banded symmetric and banded skew-symmetric matrices.
\end{itemize}
For complex matrices $X \in \mathbb{C}^{n \times n}$ we analogously consider
\begin{itemize}
\item Hermitian matrices, $X=X^*$;
\item skew-Hermitian matrices, $X=-X^*$;
\item orthogonal matrices, $X X^*= I_n$;
\item banded Hermitian and banded skew-Hermitian matrices.
\item singular values of rectangular matrices, real and complex. 
\end{itemize}
The general technique is to count the real dimension of each structured matrix manifold $\cN$, corresponding to the classes of matrices listed above, via their structured eigendecomposition. Then, we will do the same for the submanifolds of derogatory matrices of the same structure. In doing so, special attention must be paid to the redundancies in the eigendecomposition. For example, $2 \times 2$ derogatory Hermitian matrices clearly have real dimension $1$, being necessarily of the form $a I_2, a \in \R$. On the other hand, their structured decomposition is $Q (a I_2) Q^*$ for \emph{any} unitary matrix $Q$, apparently suggesting a real dimension of $5$. The clue is that we need to subtract the $4$ degrees of freedom corresponding to the redundancy in this representation. This redundancy stems from the fact that we can multiply the eigenspace by any $2 \times 2$ unitary matrix $U$ without changing the corresponding matrix.

\section{Full real matrices}
\subsection{Real symmetric matrices}

For real symmetric matrices we expose the beautiful geometric argument given by Lax~\cite{Laxlaa}, who credits Wigner and von Neumann~\cite{von1993verhalten}. The dimension of the ambient space of $n \times n$ real symmetric matrices is manifestly $\frac{n(n+1)}{2}$, while to count the dimension of the manifold $\cM$ of real symmetric matrices with at least two identical eigenvalues we proceed as follows. There are $n-1$ degrees of freedom to choose the eigenvalues. We then fix the eigenvectors corresponding to the simple eigenvalues. The first one depends on $n-1$ parameters ($n$ elements minus the constraint on unit length), the second one on $n-2$ (unit length, and orthogonality), $\dots$, the last, i.e., $(n-2)$th, one on $2$ parameters. At this point, the eigenspace corresponding to the double eigenvalue is uniquely determined. This gives a total of

$$ n-1 + \sum_{k=2}^{n-1} k = n-1 + \frac{n(n-1)}{2} -1 = \frac{n(n+1)}{2} - 2.$$
Note that we are safely neglecting the redundancies in the signs of the eigenvectors corresponding to simple eigenvalues, as these correspond to (zero-dimensional) binary choices between $-1$ and $1$. On the contrary, there is a $1$-dimensional redundancy in the choice of an orthonormal basis for the double eigenspace, which explains why everything is determined at that point.

Hence, the submanifold of real symmetric matrices with multiple eigenvalues has a real codimension $2$. For this reason, by Colollary \ref{cor:thom}, any regular curve which is a submanifold of the manifold of real symmetric matrices will generically not intersect the submanifold of symmetric matrices with double eigenvalues: one needs at least two parameters.

\begin{figure}[ht]
\centering
\includegraphics[width=70mm]{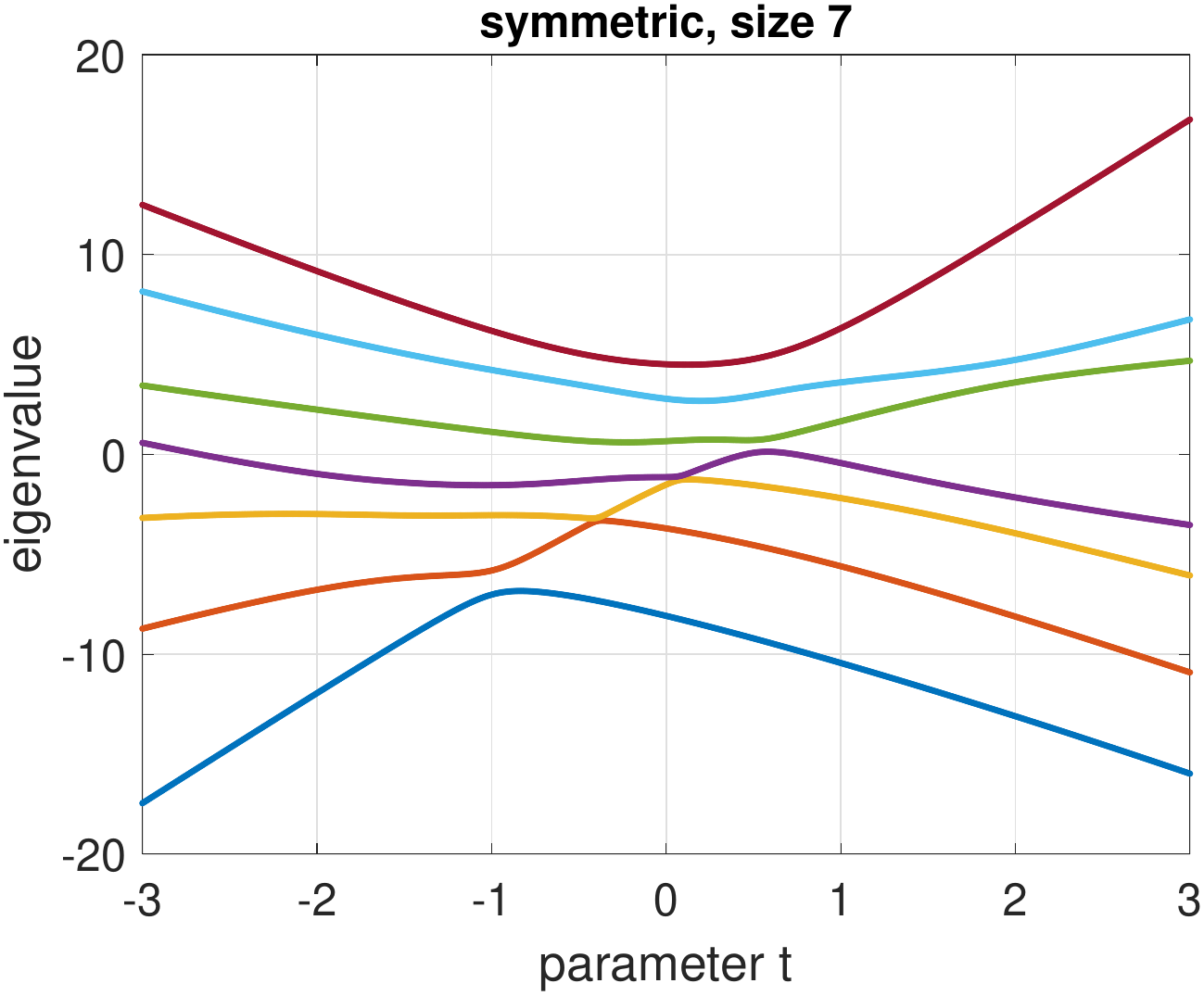}
\caption{Avoidance of the eigenvalues of a $7 \times 7$ real symmetric curve $S(t)=A+tB$ where $A,B$ are symmetric matrices. At the two points where two curves appear to cross, they become very close but do not intersect. Codimension 2. 
}
\label{fig1}\end{figure}

\subsection{Real skew-symmetric matrices}

Before generalizing Lax's argument, it is useful to recall a few properties of the spectral eigendecomposition of a generic real skew-symmetric matrix. We will then give an argument on why eigenvalue avoidance may or may not happen for this structure according to the parity of the size.

\subsubsection{Even size}\label{sec:evenskew}
Suppose that $K=-K^T \in \R^{2m \times 2m}$.
The eigenvalues are always pure imaginary and appear in complex conjugate pairs $\{ \lambda,-\lambda\}$. As a consequence, the eigenvectors also appear in complex conjugate pairs. Since the spectral theorem applies, we have an eigendecomposition of the form $K = V^* \Lambda V$ with

$$ V= \begin{bmatrix} U_1 + i U_2 & U_1 - i U_2\\
U_3 + i U_4 & U_3 - i U_4 \end{bmatrix}, \qquad \Lambda = \begin{bmatrix} i \Lambda_+ & \\
& -i \Lambda_+ \end{bmatrix} $$
where $\Lambda_+ \in \R^{m \times m}$ is real diagonal positive semidefinite. Moreover, it is easy to show that  $V$ is unitary if and only if $U$ is (real) orthogonal, where 

\begin{equation*}  
 U = \sqrt{2} \begin{bmatrix} U_1 & U_2\\
U_3 & U_4 \end{bmatrix} \in \R^{2m \times 2m}.  
\end{equation*}

Generically, we are allowed $m$ choices for $\Lambda_+$. Choosing $V$ is equivalent to choosing $U$, up to $m$ arbitrary choice of phases (one per each eigenvector in the left half of $V$). The choices of phases are more significant than the arbitrary choices of sign of a column of an orthogonal matrix --- to explain the analogy, note that the positive real semiaxis has the same dimension as the real axis, but requiring a complex number to have a fixed phase reduces the (real) dimension of the space of potential choices by one. This observation leads to the following dimensional counting argument: the first column of $V$ depends on $(2m-1)+(2m-2)-1= 4m-4$ real parameters; the second column of $V$ depends on $(2m-3)+(2m-4)-1=4m-8$ real parameters; $\dots$; the $(m-1)$th column of $V$ depends on $4$ real parameters, and one can check that at this point all the other eigenvectors are completely determined. In summary, the dimension of the ambient space of skew-symmetric matrices of even size $2m$ is

$$m + 4 \sum_{k=1}^{m-1} k = m + 2m(m-1) = m (2m- 1).$$

Of course, this could have been obtained directly, but counting the dimension via the spectral decomposition is a useful exercise as we are now 
going to modify the argument slightly to estimate the size of the submanifold of skew-symmetric matrices with even size and at least one multiple eigenvalue.

What is the dimension of the submanifold $\cM$ of real skew-symmetric matrices of even size $2m$ with at least one multiple eigenvalue? The answer depends on where the collision happens. If the double eigenvalue is at $0$, we have one degree of freedom less in determining $\Lambda_+$, but the argument above is otherwise unchanged, so that $\cM$ has real codimension $1$. An alternative proof, less close in spirit to Lax's and Keller's argument, is the following: the set of skew-symmetric singular matrices has obviously codimension $1$, as it is the locus of the equation $\det K = 0$. Because of the spectral symmetry of skew-symmetric real matrices, however, the latter set is equal to the set of skew-symmetric matrices of even size that have a $0$ eigenvalue of multiplicity at least $2$.

If the collision happens at an eigenvalue $i \mu \neq 0$ (and hence, also at $-i \mu$), which of course is possible only if $n = 2m  \geq 4$, then we have $m-1$ degrees of freedom in choosing $\Lambda_+$; moreover, once the $(m-2)th$ eigenvector of $V$ has been fixed, we are left with a subspace of dimension $4$ that is spanned by the eigensubspaces corresponding to $\pm i \mu$. We have some more freedom in selecting a $2$-dimensional complex eigenspace for $ i \mu$, with the constraint of having the special structure coming from the fact that it is an eigenspace of a \emph{real} skew-symmetric matrix. The argument to count the real degrees of freedom is similar to the one above: this time, we count the degrees of freedom to choose $4$ columns of a real orthogonal matrix $U$, from which we subtract the ``illusory'' degrees of freedom corresponding to the arbitrary choise of postmultiplying by any $2 \times 2$ complex unitary matrix. We can imagine to have projected on the orthogonal complement of the previously selected eigenvectors, and hence the matrix sizes are fixed to be $4$: this gives $\frac{4 \cdot 3}{2} - 2^2 = 2$ more degrees of freedom.
Hence, the codimension of $\cM$ is $1+4-2=3$.

\begin{figure}[ht]
\centering
\includegraphics[width=70mm]{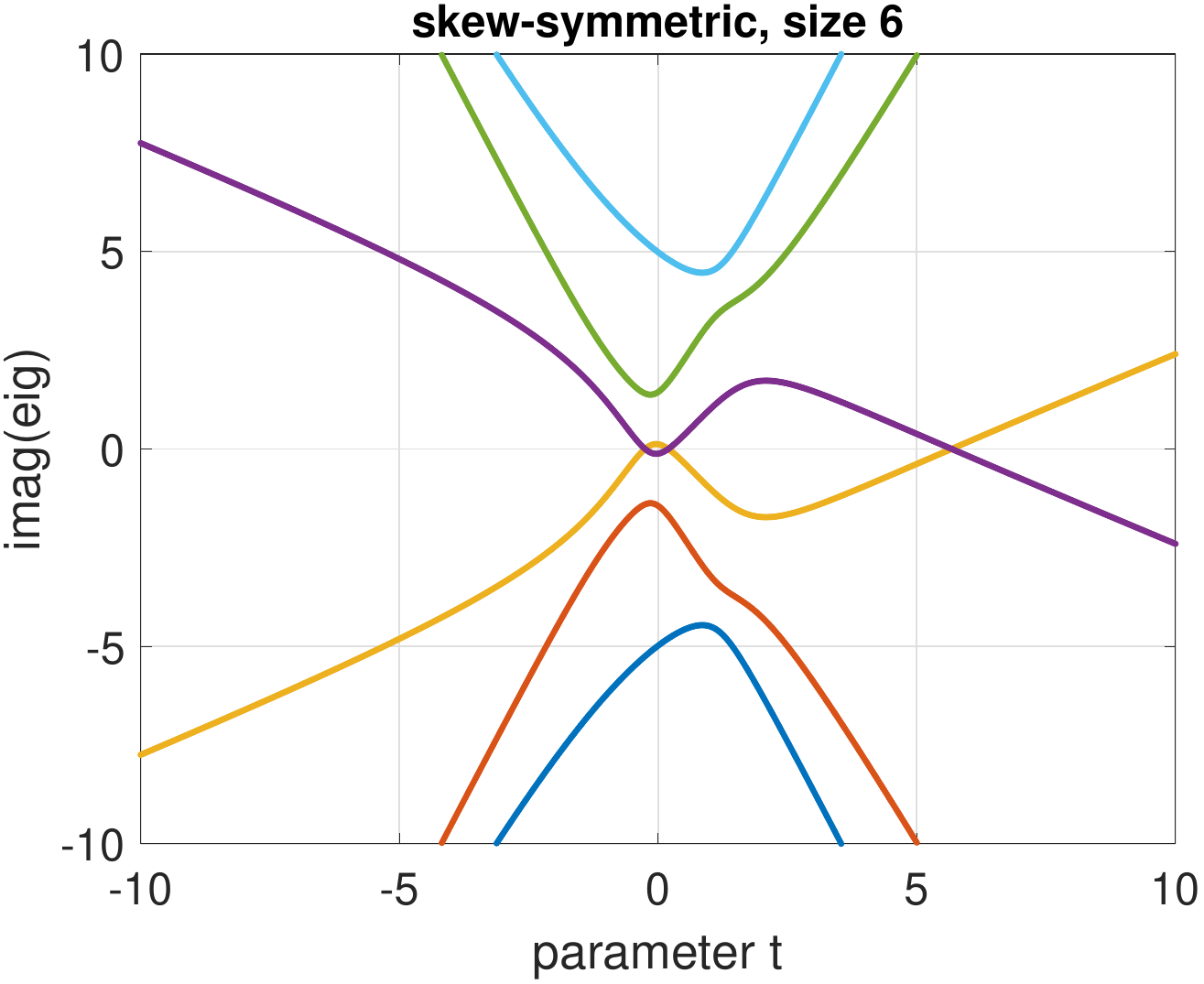}
\caption{Eigenvalues (imaginary parts) of a $6 \times 6$ real skew-symmetric curve $K(t)=A+tB$. Avoidance is seen everywhere, except at $\lambda=0$, where a crossing does happen. Codimension 1 (3 for $\lambda \neq 0$).}
\label{fig2}\end{figure}

\subsubsection{Odd size}

Let now $K=-K^T \in \R^{(2m+1) \times (2m+1)}$. For odd size, $K$ must always be singular, and the zero eigenvalue is associated with a real eigenvector. Hence, the eigenvalues and eigenvectors have the form, respectively,

$$ \Lambda=\begin{bmatrix}
i \Lambda_+ & & \\
& 0 & \\
& & -i \Lambda_+
\end{bmatrix}, \qquad V=\begin{bmatrix}
U_1+i U_2 & U_0 & U_1 - i U_2 \end{bmatrix} $$
where $\Lambda_+$ is real diagonal positive semidefinite, $U_0 \in \R^{2m+1}$, and the matrix

$$ U=\begin{bmatrix} \sqrt{2} U_1 & U_0 & \sqrt{2} U_2 \end{bmatrix}$$
is real orthogonal.

Generically, we have $m$ real degrees of freedom for $\Lambda_+$, and again choosing $V$ is equivalent to choosing $U$, up to $m$ degrees of freedom lost for arbitrary choices of phases (note that we already implicitly made a phase choice for $U_0$ by claiming that it is real). $U_0$ depends on $2m$ real parameters, and then the other eigenvectors of $V$ depend, as before, on $4m-4, 4m-8, \dots, 4$ parameters. After having fixed $U_0$ and the first $m-1$ columns of $V$, all the other entries are determined. The dimension of the ambient space is therefore, as expected,

$$ m + 2 m + 4 \sum_{k=1}^{m=1} k = m (2 m + 1).$$

The analysis for a collision at a nonzero eigenvalue is similar to before. The real degrees of freedom for the eigenvalues become $1$ less, and the real degrees of freedom for the eigenvectors become $4-2=2$ less. Therefore, the codimension is $3$ and we do see eigenvalue avoidance.

When a collision happens at zero the number of columns of $U_1,U_2,\Lambda_+$ shrinks from $m$ to $m-1$ while the number of columns of the zero diagonal block in $\Lambda$ and $U_0$ grows from $1$ to $3$. Hence, there are $m-1$ degrees of freedom for the eigenvalues. The first eigenvector of $V$ depends on $(2m)+(2m-1)-1=4m-2$ real parameters, the second on $(2m-2)+(2m-3)-1=4m-6$ real parameters, and so on until the $(m-1)th$ eigenvector of $V$ that depends on $6$ real parameters. At this point, every column of $V$ is determined except for those of $U_0$, but the eigenspace of $0$ is specified, so the dimension counting gives

$$ m-1 + \sum_{k=1}^{m-1} (4k+2)= (m-1)(2m+3).$$
Thus the codimension is $m(2m-1)-(m-1)(2m+3)=3$.

We summarize our findings in the following theorem.

\begin{theorem}\label{thm:skewsym}
The codimension of derogatory real $n \times n$ skew-symmetric matrices, embedded in the manifold of real $n \times n$ skew-symmetric matrices, is $1$ if $n$ is even and $3$ if $n$ is odd.

For $n \geq 3$, the codimension of real $n \times n$ skew-symmetric matrices that have at least a double nonzero eigenvalue, embedded in the manifold of real $n \times n$ skew-symmetric matrices, is $3$.
\end{theorem}

To clarify, the codimension 1 in the even case comes entirely from the possibility of the multiple eigenvalue at 0. For $n \leq 2$, there cannot be a double nonzero eigenvalue. 

\begin{figure}[ht]
\centering
\includegraphics[width=70mm]{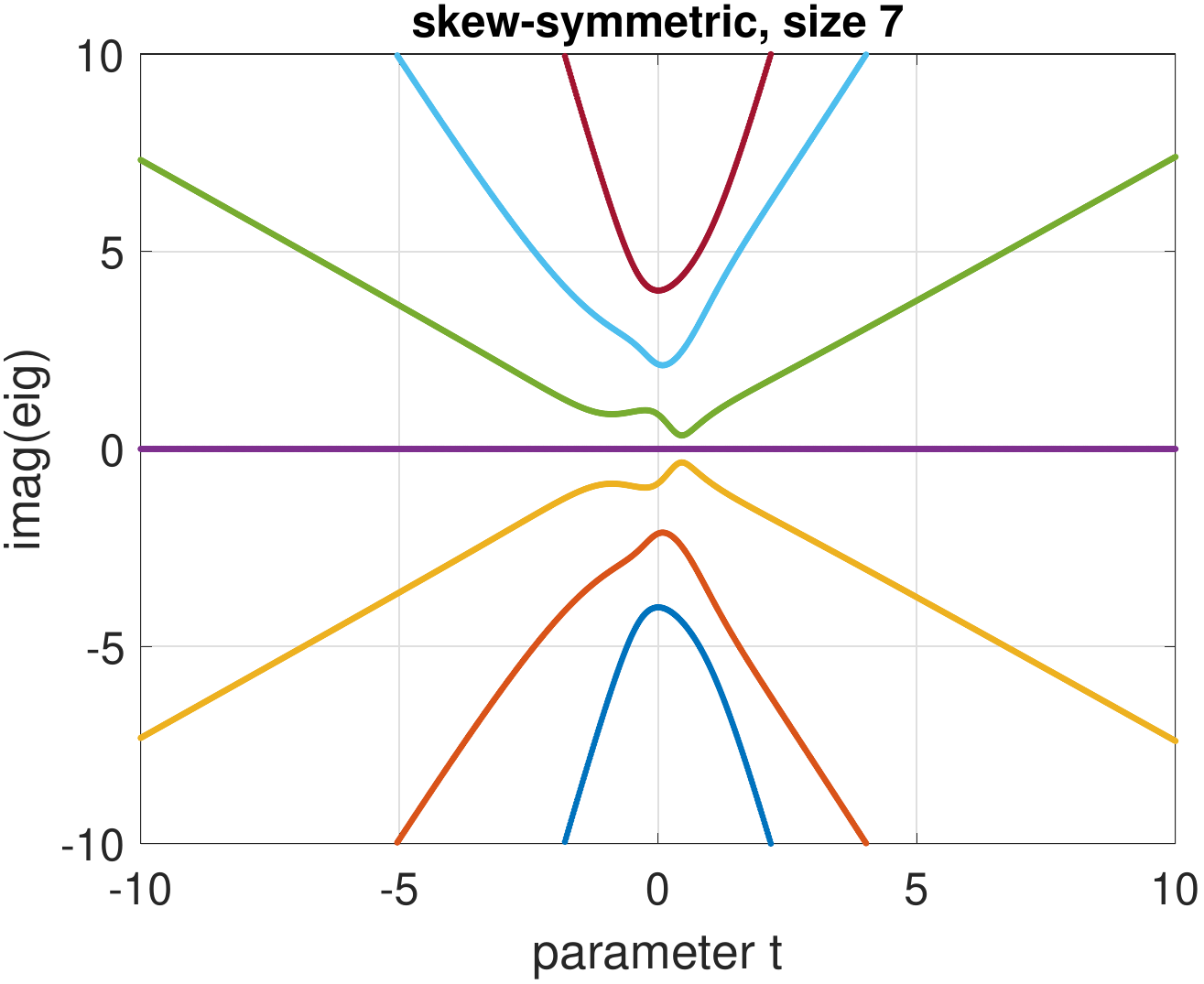}
\caption{Eigenvalues of a $7 \times 7$ real skew-symmetric curve $K(t)=A+tB$. Unlike for even sizes, for odd sizes avoidance is expected (and seen) also at zero, i.e., the nonzero eigenvalues are repelled by the constantly zero unpaired eigenvalue. Codimension 3. }
\label{fig3}\end{figure}

\subsection{Real orthogonal matrices}

Let $X \in \R^{n \times n}$ be orthogonal. If $n$ is odd or $\det(X)=1$ (or both), then there exists a skew-symmetric $K \in \R^{n \times n}$ such that $\det(X)X= e^K$.  To see this, note that in this setting $\det(\det(X)X)=1$, and hence there is a unitary matrix $U$ such that $\det(X)X=UDU^*$ where $D$ is diagonal and for some $\theta_1,\dots,\theta_m \in ]-\pi,\pi]$ it holds
\[  D = \begin{cases}
\mathrm{diag}(e^{i \theta_1},e^{-i\theta_1},\dots,e^{i \theta_m},e^{-i \theta_m},1) \ &\mathrm{if} \ n=2m+1 \ \mathrm{is} \ \mathrm{odd};\\
\mathrm{diag}(e^{i \theta_1},e^{-i\theta_1},\dots,e^{i \theta_m},e^{-i \theta_m}) \ &\mathrm{if} \ n=2m \ \mathrm{is} \ \mathrm{even}.
\end{cases}  \]
Hence, we can costruct
\[  K = \begin{cases}
U\mathrm{diag}(i \theta_1,-i\theta_1,\dots,i \theta_m,-i \theta_m,0)U^* \ &\mathrm{if} \ n=2m+1 \ \mathrm{is} \ \mathrm{odd};\\
U\mathrm{diag}(i \theta_1,-i\theta_1,\dots,i \theta_m,-i \theta_m)U^* \ &\mathrm{if} \ n=2m \ \mathrm{is} \ \mathrm{even}.
\end{cases}  \]
For this reason, under the above assumptions, the eigenvector structure of real orthogonal matrices is precisely the same as that of real skew-symmetric matrices. However, this is not so simple when $n$ is even and $\det(X)=-1$, which requires a separate analysis. Moreover, the spectrum of orthogonal matrices lie on the unit circle, which (unlike the imaginary axis) intersect the real axis on two, and not just one, special points: $+1$ and $-1$. Finally, geometrically real orthogonal matrices are a manifold with two connected components, $\det X=1$ and $\det X=-1$. Coherently with all these observations, it turns out that their behaviour with respect to eigenvalue avoidance depends both on the parity of the size and on the parity of the determinant.

\subsubsection{Even size}

Consider first matrices belonging to the connected component of real orthogonal matrices having determinant $1$. The analysis is fairly similar to the case of skew-symmetric matrices of even size, except that there are two exceptional points. The conclusion is that collisions at eigenvalues $\neq \pm 1$ happen in a submanifold of codimension $3$, while collisions at $\pm 1$ happen in a submanifold of codimension $1$, and hence there is no eigenvalue avoidance for curves within orthogonal matrices of even size. This is similar to what happens for skew-symmetric matrices of even size, but collisions are expected at two special points rather than just at one special point.

If the determinant is $-1$, the situation is different in the sense that one eigenvalue is always $1$ and another eigenvalue is always $-1$. The spectral theorem yields an eigendecomposition of the form $X=V^* \Lambda V$ with
$$ \Lambda = \begin{bmatrix}
e^{i \Lambda_+} & & & \\
& 1 &  & \\
& & -1 & \\
& & & e^{-i \Lambda_+}
\end{bmatrix}, \qquad  V=\begin{bmatrix}
U_2 + i U_3 & U_0 & U_1 & U_2-iU_3
\end{bmatrix}  $$
where $\Lambda_+$ is positive semidefinite with eigenvalues in $[0,\pi]$, and the matrix $U=\begin{bmatrix}
\sqrt{2} U_2 & U_0 & U_1 & \sqrt{2} U_3
\end{bmatrix}$ is real orthogonal.
The dimension count yields, similarly to the analysis for skew-symmetric matrices and for $n=2m$,
$$(m-1)+(2m-1)+(2m-2)+\sum_{k=1}^{m-1}(4m-4(k+1)) = 2m^2-m$$
as expected.

Nothing changes, in essence, in the analysis for collisions at eigenvalues $\neq \pm 1$, which still happens in a submanifold of codimension $3$. Collisions at $\pm 1$ can happen only if the size is at least $4$. Noting that the number of columns of $U_2,U_3$ shrink by one while the number of columns of either $U_0$ or $U_1$ (but not both) increase by two, we are led to an argument which is similar to that for the collision at $0$ for \emph{odd}-sized skew symmetric matrices. The conclusion is that we see a loss of $1$ degree of freedom for eigenvalues and $2$ for eigenvectors. Therefore collisions at $\pm 1$ happen in a submanifold of codimension $3$. We expect repulsion of the non-constant eigenvalue functions from the eigenvalues constantly equal to $\pm 1$.

\begin{figure}[ht]
  \begin{minipage}[t]{0.5\hsize}
\includegraphics[width=.9\textwidth]{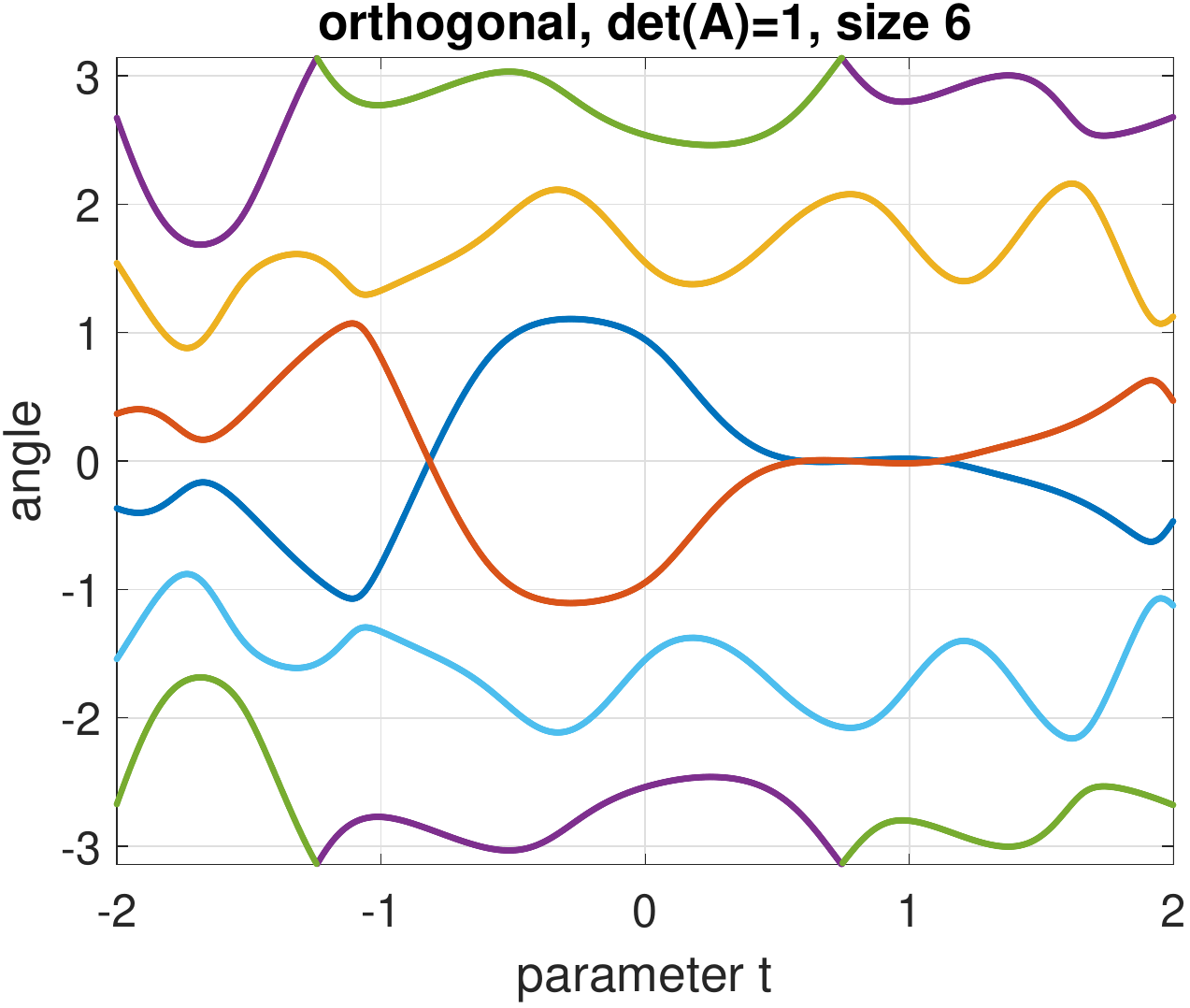}
  \end{minipage}   
  \begin{minipage}[t]{0.5\hsize}
\includegraphics[width=.9\textwidth]{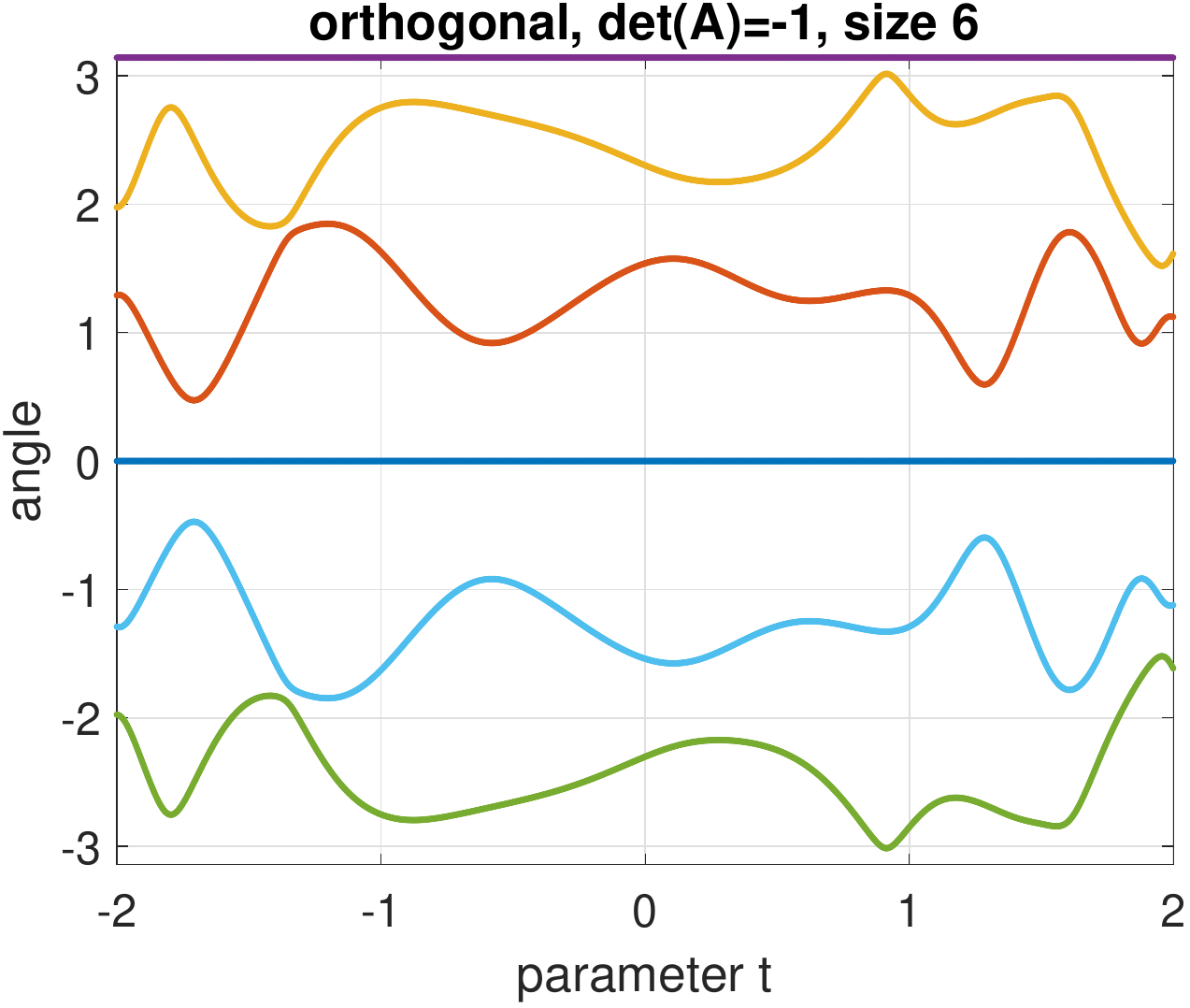}
  \end{minipage}
\caption{
Left: Eigenvalues (their arguments $\in(-\pi,\pi]$) of a $6 \times 6$ real orthogonal curve $Q(t)$ with $\det Q(t) \equiv 1$, obtained by the orthogonal factor in the polar decomposition of $Q(0)(I+tS)$, where $S$ is a skew-symmetric matrix (see Section~\ref{sec:unit} for the construction). 
Collisions happen at the eigenvalues $\pm 1$ (angles $0$ and $\pi$). 
 Codimension 1 ($\lambda=\pm 1$), 3 (elsewhere). 
Right: Eigenvalue avoidance for a $6 \times 6$ real orthogonal curve $Q(t)$ with $\det Q(t) \equiv -1$. The experiment confirms the prediction of a repulsion from the eigenvalues constantly equal to $\pm 1$ (angles $0$ and $\pi$).
Codimension 3. 
}
\label{fig6}\end{figure}

\subsubsection{Odd size}

Note that the manifold of real orthogonal matrices has two connected component, according to the value of the determinant, $+1$ or $-1$. Any odd ($2m+1$) size orthogonal matrix with determinant $\delta$ has at least one eigenvalue at $\delta$. We partially replicate the analysis for skew-symmetric cases of odd size: collisions at eigenvalues not equal to $\pm \delta$ happen in a submanifold of codimension $3$, and triple collisions at $\delta$ happen in a submanifold also of codimension $3$. There is a third case that was not present for skew-symmetric matrices though: one can have a collision at $-\delta$. There is a degree of freedom lost for eigenvalues, while for eigenvectors we can see that, as for real skew-symmetric matrices, once we have chosen $m$ eigenvectors corresponding to eigenvalues generically not equal to $\pm \delta$, everything else is determined already: therefore, there is no loss of degrees of freedom, and the codimension is $1$. We conclude that, unlike for skew-symmetric matrices of odd size, we do not expect eigenvalue avoidance. Also, unlike orthogonal matrices of even size, collisions are expected to happen at only one special point, equal to minus the determinant of the matrix.

\begin{figure}[ht]
  \begin{minipage}[t]{0.5\hsize}
\includegraphics[width=.9\textwidth]{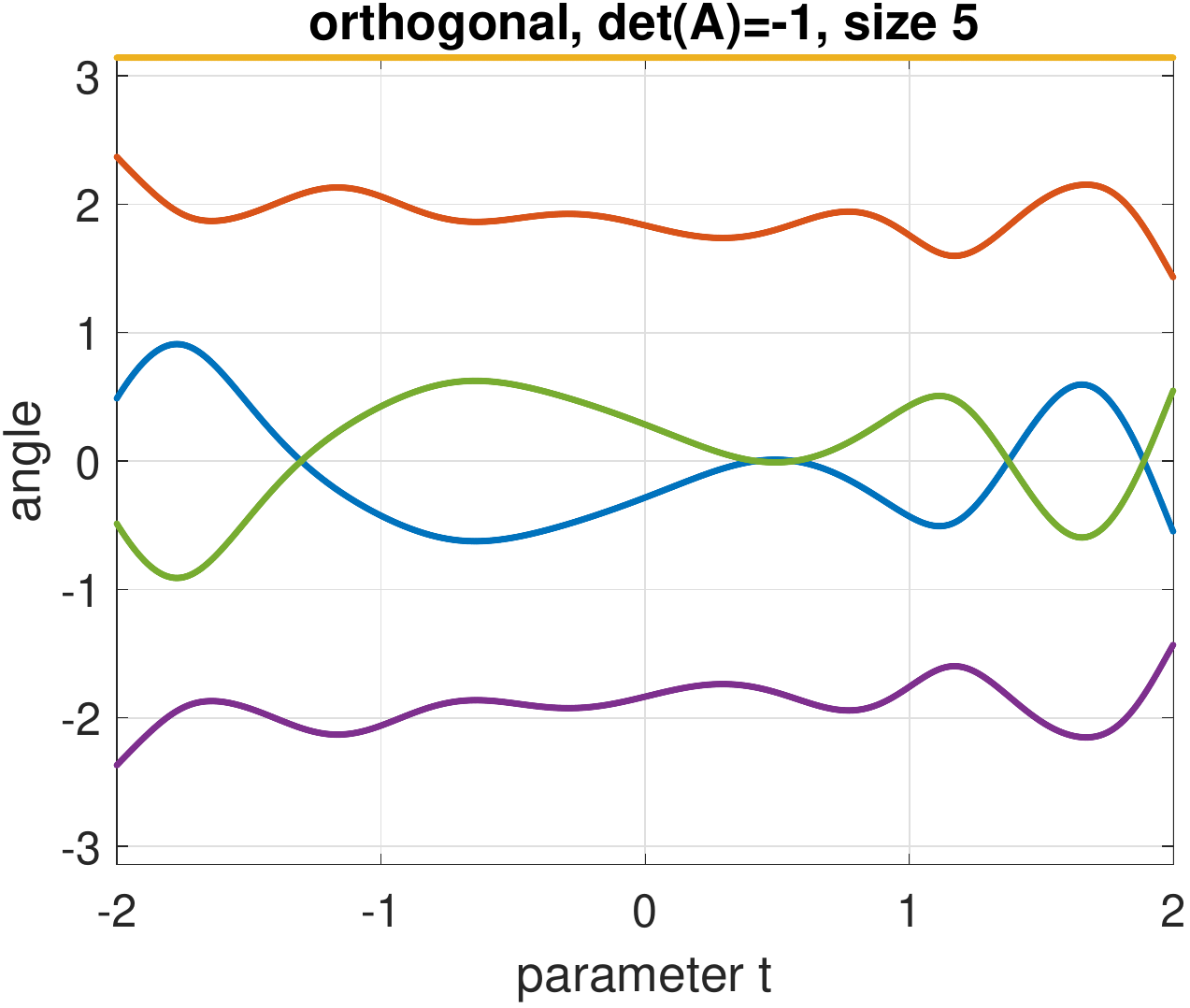}
  \end{minipage}   
  \begin{minipage}[t]{0.5\hsize}
\includegraphics[width=.9\textwidth]{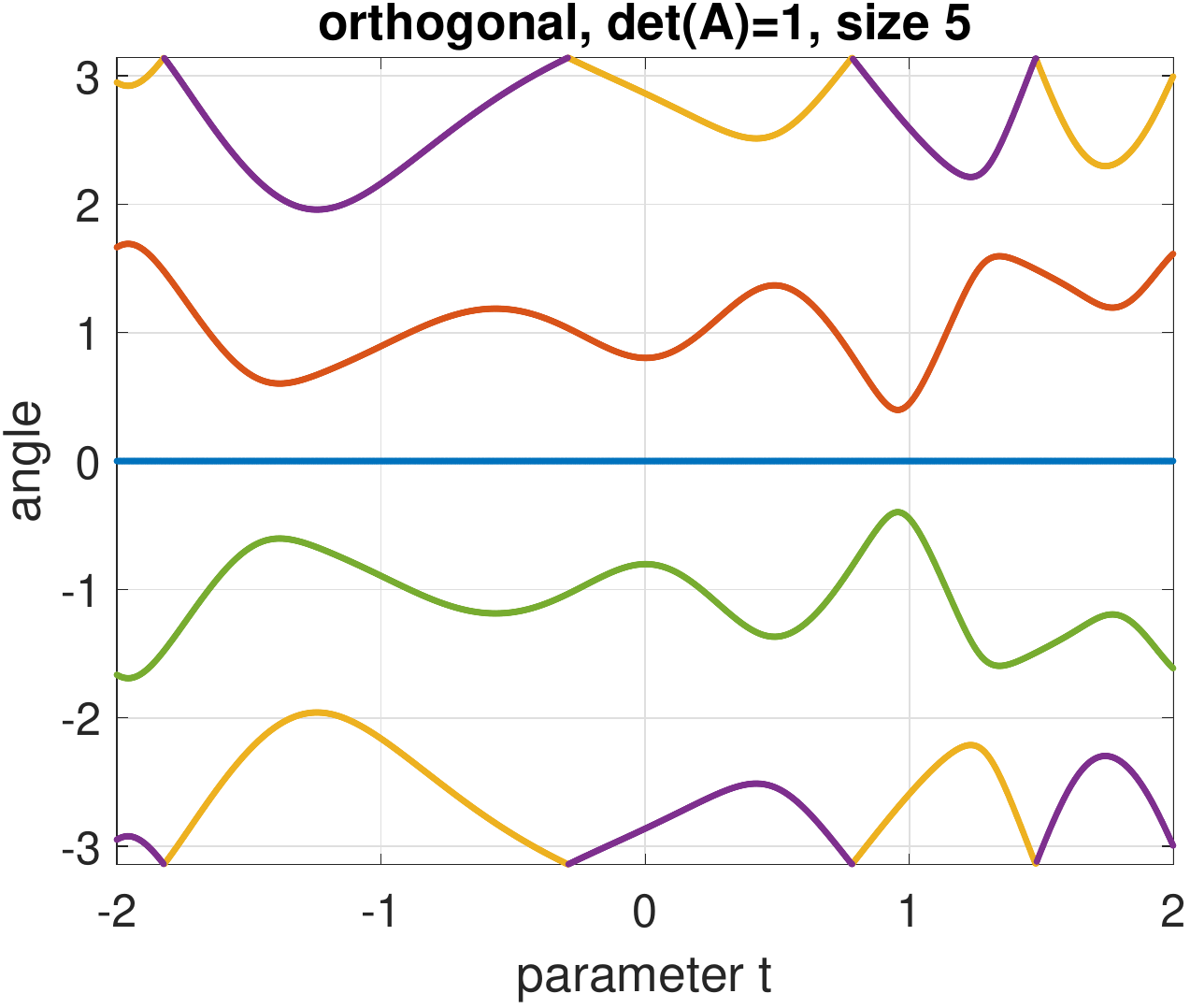}
  \end{minipage}
\caption{Eigenvalues of a $5 \times 5$ real orthogonal curve $Q(t)$ with $\det Q(t) \equiv 1$, generated as in Figure~\ref{fig6}.
As predicted by the theory, we see avoidance everywhere except at the eigenvalue $1$ (angle $0$) when $\det Q(t)=-1$ (left), and at  $-1$ (angle $\pi$) when $\det Q(t)=1$ (right). 
Codimension is 3 except when $\det Q(t)=-1$ and at $\lambda=1$.
}
\label{fig7}\end{figure}

To summarize our analysis of orthogonal matrices:

\begin{theorem}\label{thm:orth}
The codimension of derogatory real $n \times n$ orthogonal matrices with determinant $\delta$, embedded in the manifold of real $n \times n$ orthognal matrices with determinant $\delta$, is
\begin{itemize}
\item $1$ if $n$ is even and $\delta=1$;
\item $3$ if $n$ is even and $\delta=-1$;
\item $1$ if $n$ is odd.
\end{itemize}

Moreover, if $n \geq 4$ is even, the codimension of real $n \times n$ orthogonal matrices with determinant $\delta$ with at least one double nonreal eigenvalue is $3$; and if $n$ is odd, the codimension of real $n \times n$ orthogonal with determinant $\delta$ with at least one double eigenvalue not equal to $-\delta$ is $3$.
\end{theorem}

\section{Full complex matrices}

\subsection{Hermitian matrices}

\begin{theorem}\label{thm:hermitian}
The real submanifold of derogatory Hermitian matrices, when embedded in the real manifold of Hermitian matrices, has real codimension $3$.
\end{theorem}

The statement of Theorem \ref{thm:hermitian} is given by Lax \citep{Laxlaa} but its proof is left as an exercise for the reader~\citep{Laxlaa}. 
A solution is given by Keller in~\cite{keller2008}.
We briefly give here our own solution to that exercise.  Of course, the real dimension of the real manifold of Hermitian matrices of size $n$ is $n^2$. Let us now count the dimension of the submanifold $\cM$ of Hermitian matrices having at least one multiple eigenvalue. Of course the degrees of freedom for the eigenvalues are 
 $n-1$. Unlike for real vectors, the orthogonality conditions for complex vectors give $2$ real constraints each; the unit length conditions is still equivalent to $1$ real constraint. So we have $2n-2$ real degrees of freedom for the first eigevector ($2n$ elements, subtractring one for the unit length constraint and one for arbitrarily fixing the phase), $2n-4$ for the second, $\dots$, $4$ for the $(n-2)$th, and at this point the eigenspace of the multiple eigenvalue is uniquely determined. Hence,

$$\dim M = n-1 + 2 \sum_{k=1}^{n-2} (n-k) = n^2-3.$$

\subsection{Skew-Hermitian matrices}

\begin{theorem}\label{thm:skewHerm}
The real submanifold of derogatory skew-Hermitian matrices, when embedded in the real manifold of skew-Hermitian matrices, has real codimension $3$.
\end{theorem}

\begin{proof}
There is a bijection between skew-Hermitian and Hermitian matrices given by $M \mapsto i M$. Therefore, the situation is precisely the same as in the previous subsection, except of course that the eigenvalues lie on the imaginary axis rather than on the real axis, and there is nothing else to say.
\end{proof}

\begin{figure}[ht]
  \begin{minipage}[t]{0.5\hsize}
\includegraphics[width=.9\textwidth]{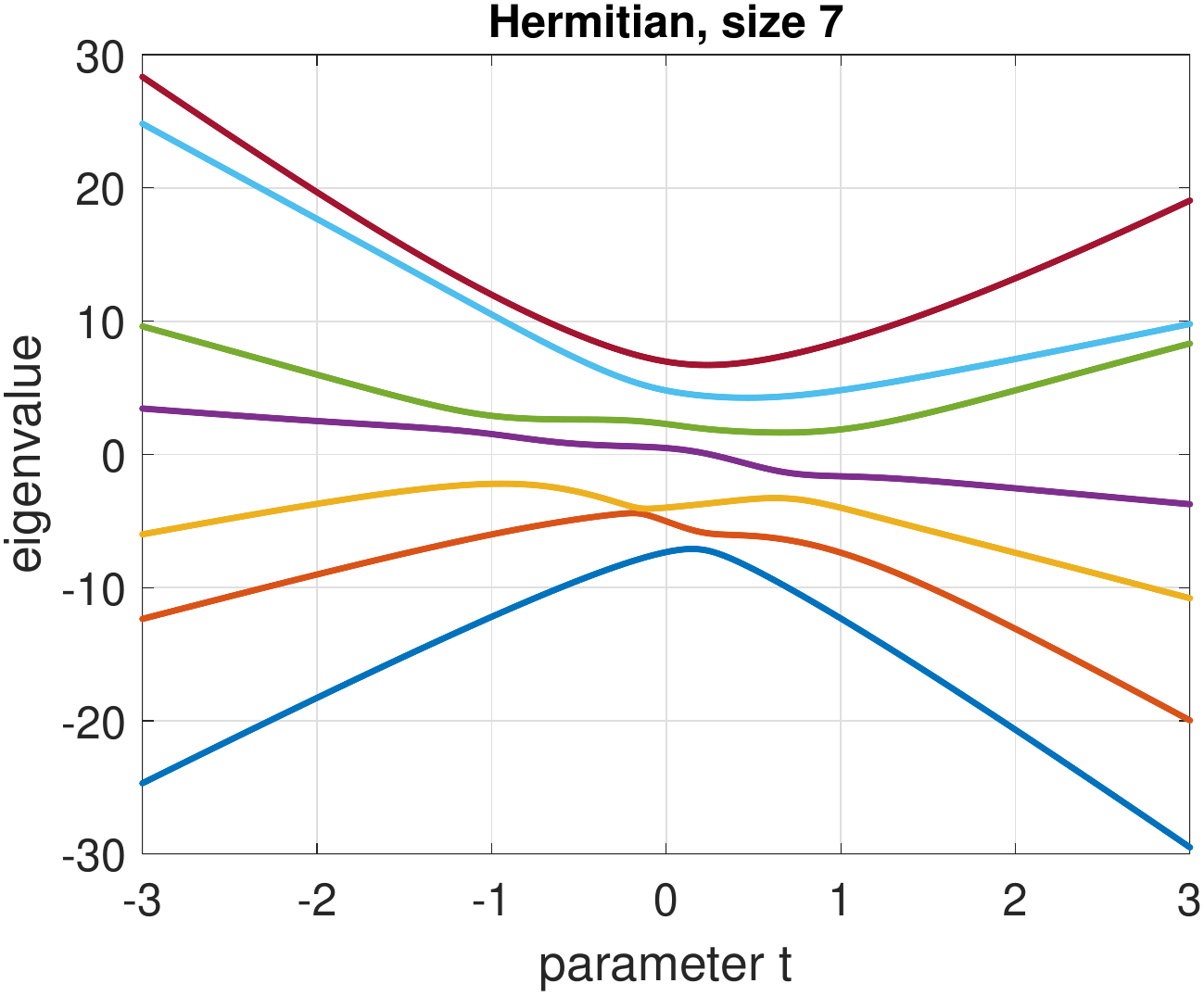}
  \end{minipage}   
  \begin{minipage}[t]{0.5\hsize}
\includegraphics[width=.9\textwidth]{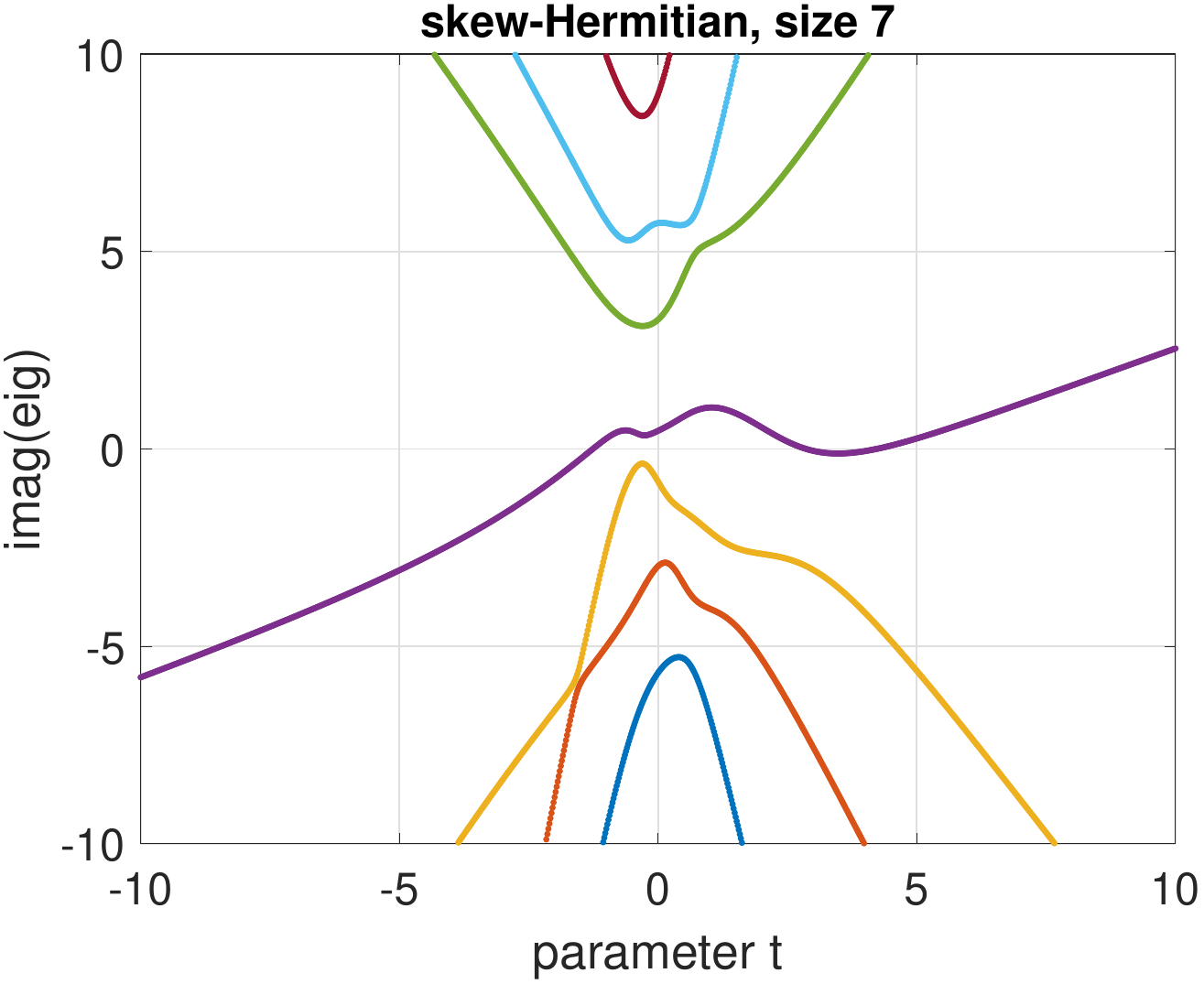}
  \end{minipage}
\caption{Eigenvalue avoidance for a $7 \times 7$ 
Hermitian (left) and skew-Hermitian (right) curves $W(t)=A+tB$. Both have codimension 3.}
\label{fig5}\end{figure}

\subsection{Unitary matrices}\label{sec:unit}

If $U \in \C^{n \times n}$ is unitary then $U=e^{i H}$ for some Hermitian $H \in \C^{n \times n}$. Again, the dimension counting is perfectly similar to Hermitian matrices, and any collision of the eigenvalue functions happens in a subspace of codimension $3$.

It is worth noting that if one has a curve $H(t)$ of Hermitian matrices, the curve of unitary matrices $U(t)=e^{i H(t)}$ is likely to have collisions. This is because it suffices that $H(t)$ has eigenvalues that differ by $2 k \pi$ (which happens in a submanifold of codimension $1$) for $U(t)$ to have multiple eigenvalues. This is, however, an artificial consequence of the particular construction of this path 
(note that the exponential map is not injective and hence that this curve may not be a submanifold, which is a situation beyond the applicability of Corollary \ref{cor:thom}) 
as can be verified by building the path differently, for example, $U(t) = (I - i H(t))(I + i H(t))^{-1}$. Observe that this rational map is injective, but not surjective: it maps Hermitian matrices onto $\{UU^*=I, \det  (U+I) \neq 0\}$, that is, no eigenvalue can be $-1$.
Alternatively, to avoid the artificial eigenvalue avoidance at $-1$, one can obtain a path via the unitary polar factor $U_p(t)$ of the matrix $Q(0)(I+tS)$, where $Q(0)$ is unitary, $S\neq 0$ is skew-Hermitian. We use this in our experiments. 
Note that $Q(0)(I+tS)$ lies in the tangent space of unitary matrices~\cite{edel98}; and 
for all $t$ we have
$\det(I+tS)\neq 0$ (by the Fan-Hoffman theorem~\cite[Prop.~III.5.1]{bhat:96}, $\sigma_i(I+tS)\geq 1$ for all $t\in\mathbb{R}$, so $U_p(t)$ is uniquely defined and continuous in $t$. One can also see that $t\mapsto U_p(t)$ is injective.

In summary, we have the following. 

\begin{theorem}
The real submanifold of derogatory unitary matrices, when embedded in the real manifold of unitary matrices, has real codimension $3$.
\end{theorem}

\begin{figure}[ht]
\centering
\includegraphics[width=70mm]{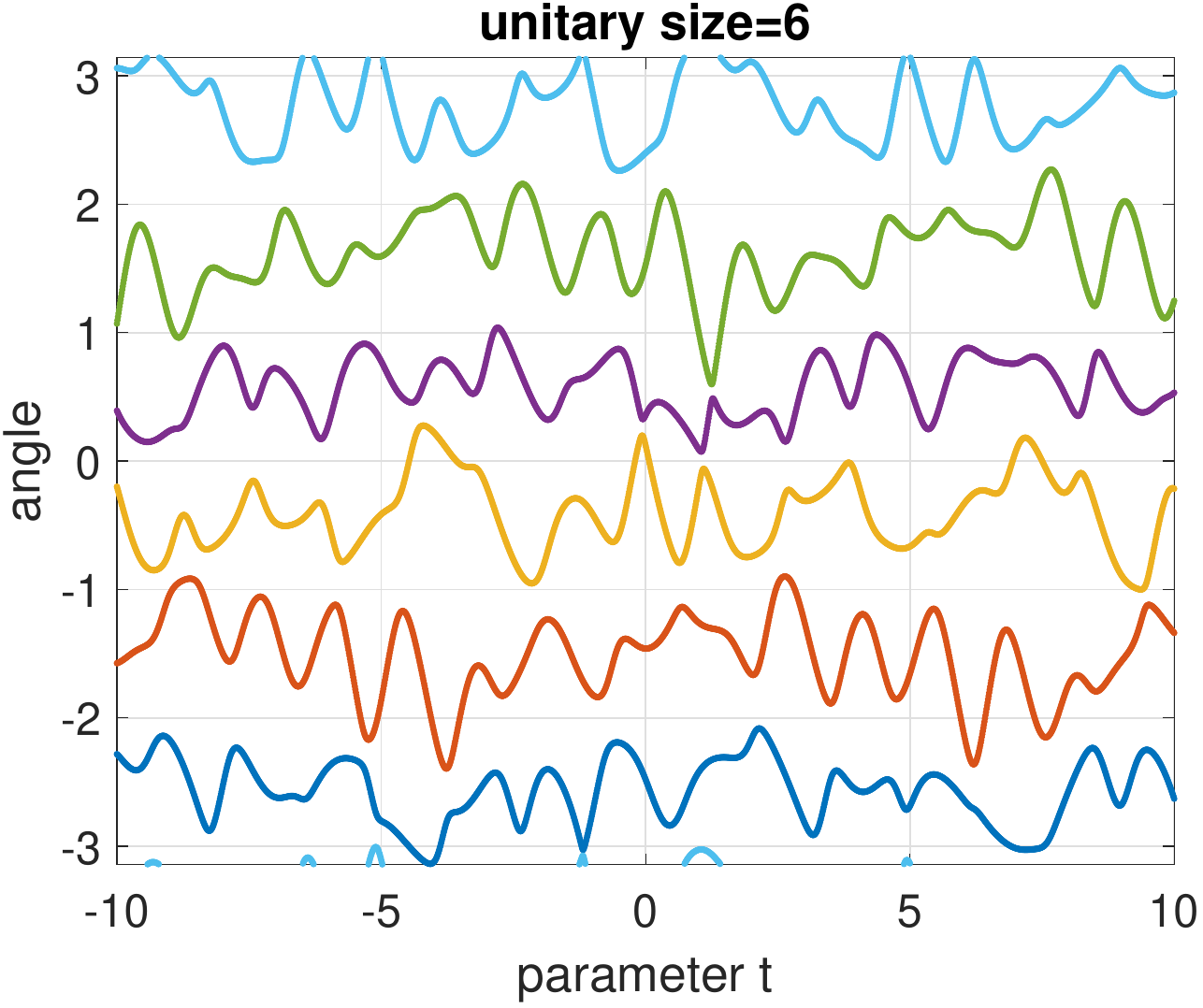}
\caption{Eigenvalue avoidance for a $6 \times 6$ unitary curve $U(t)$. Codimension 3. }
\label{fig5bis}\end{figure}

\section{Banded matrices}

Recall that the bandwidth of a matrix with elements $a_{ij}$ is defined as $ \max |j-i| : a_{ij} \neq 0.$ For example, nonzero diagonal matrices have bandwidth $0$, and nondiagonal tridiagonal matrices have bandwidth $1$.

We will analyze here the case of structured matrices with bandwidth $k \leq n-1$, where $k=0$ corresponds to the diagonal case (which is exceptional) and $k=1,n-1$ correspond to, respectively, tridiagonal and full matrices. Banded matrices are relevant in practice, as they often arise from applications. Moreover, the standard algorithm for the symmetric (or Hermitian) eigenvalue problem is to first reduce the matrix to tridiagonal (or sometimes banded) form~\cite[Ch.~8]{golubbook4th}.

One may be tempted to conjecture that, for example, the classical result for derogatory real symmetric matrices implies that the codimension of derogatory real banded symmetric matrices is also $2$. However, this ``codimension invariance'' happens only conditionally when we intersect a manifold and its submanifold with a third manifold, and the condition has a distinctive geometric nature. What can be said is the following statement \cite[Thm.~6.30 and related discussion]{lee}.

\begin{theorem}\label{thm:thom}
Suppose that $\cM_0, \cM_1, \cN_0, \cN_1$ are smooth manifolds satisfying the following: (1) $\cM_0 \subseteq \cN_0$ has codimension $k$ in $\cN_0$ and (2) $\cM_1 = \cM_0 \cap \cN_1$. If $\cM_0 \pitchfork \cN_1$, then $\cM_1 \subseteq \cN_1$ has codimension $k$ in $\cN_1$.
\end{theorem}

In Theorem \ref{thm:thom}, $\pitchfork$ means ``intersect transversally''. In the case of our interest, we would set (for example) $\cN_0 \rightarrow$ real symmetric matrices, $\cM_0 \rightarrow$ derogatory real symmetric matrices, and $\cN_1 \rightarrow$ real symmetric matrices of bandwidth $k$. However, it is not at all clear \emph{a priori} whether $\cN_1$ and $\cM_0$ intersect transversally! Indeed, at least in the case $k=0$ (diagonal matrices), the codimension of derogatory real diagonal matrices is clearly $1$, and hence (by Theorem~\ref{thm:thom}), the intersection cannot be transversal.

For the same reason, the study of other banded structured (e.g. Hermitian, skew-symmetric, skew-Hermitian, $\ldots$) matrices will also need a case-by-case analysis, and cannot be simply inferred by the result for full matrices of the same structure.

We note that, for real symmetric and Hermitian matrices, this statement appeared in \cite{DPP18}, but not many details are given. For this reason, we feel that it is appropriate to give a detailed proof here.

\subsection{Symmetric banded matrices}

As anticipated, for this structure we have the following result.

\begin{theorem}
The codimension of derogatory banded real symmetric matrices of bandwidth $k > 0$, embedded in the manifold of real symmetric matrices of bandwidth $k$, is $2$.
\end{theorem}

\begin{proof}
If $A \in \R^{n \times n}$ is symmetric and banded of bandwidth $k>0$, then its real degrees of freedom are manifestly $$n+(n-1)+\dots+(n-k) = \frac{1}{2}(k+1)(2n-k).$$

We now need to count the dimensions of \emph{derogatory} real symmetric matrices of bandwidth $k$, and we will do this by analyzing the freedom in choosing their eigenvalues and eigenvectors, as before. 
Although the proof is essentially the same for all $k > 0$, for the sake of concreteness we first analyze $k=1$ (tridiagonal matrices) in full detail. Then, we will show how to extend the argument to any $k>1$.

We first focus on the tridiagonal case. Let us first re-count dimensions of tridiagonal real symmetric matrices via their eigendecomposition. Their eigenvalues have $n$ degrees of freedom. The freedom in the eigenvectors is exhausted by choosing \emph{one} nonzero vector in $\R^n$ modulo scaling by a scalar. To justify this claim, let $V$ be the left eigenvector matrix of $A$ such that $DV = VA$; throughout this paragraph, $V$ is orthogonal. Since $A$ is tridiagonal, this precisely represents a Krylov recurrence. In MATLAB-style notation, wherein $V(:,1:j)$ denotes the first $j$ columns of $V$,
\begin{equation}  \label{eq:kryl}
\mbox{span}(V(:,1:j)) = \mbox{span}(v_1,Dv_1,\ldots,D^{j-1}v_1),   
\end{equation}
 where $v_1$ is the first column of $V$, which is the vector of the first elements of the 
eigenvectors of $A$. Equation~\eqref{eq:kryl} indicates that once $v_1$ is chosen, the entire matrix $V$ is uniquely determined. 
\begin{remark}Strictly speaking, this process determines $A$ up to the freedom of making a choice of signs for $v_2,\dots,v_n$ (which corresponds to picking the signs of the off-diagonal elements of $A$), and up to neglecting the redundancy in the choice of the signs for the components of $v_1$, e.g., $v_1$ can be taken nonnegative without loss of generality. However, this subtlety does not contribute to the dimensional counting, and will therefore be neglected in the rest of the proof.\end{remark}
(An exception to to the argument above is if the entire Krylov subspace, when $j=n$, is not $\R^n$; this happens only with degenerate cases, i.e., in a Zariski closed proper subset, and therefore it does not contribute to the dimensionality count.) 
Choosing $v_1$ has $n-1$ degrees of freedom; hence the overall dimension is $2n-1$, as expected. 

We now modify the counting above to take care of derogatory tridiagonal matrices: we will still neglect signs and degenerate cases. The degrees of freedom for the eigenvalues are now $n-1$, since at least two of them must be equal. As before, we can pick the first components of the eigenvector ($n$ parameters), from which we subtract a normalization constraint and an arbitrary choice of phase when picking a basis within the eigenspace of the double eigenvalue, i.e., when fixing the last two elements of this vector (for example, we can make the ``canonical choice'' to rotate this basis to make sure that the last element of this vector is $0$). Everything else, up to some choices of signs, is then determined by the Krylov construction above. Hence, this gives $n-2$ degrees of freedom for the eigenvectors. As a result, the dimension of the manifold of derogatory real tridiagonal symmetric matrices is 
$$ (n-1)+(n-2)=2n-3,$$
and therefore its codimension when embedded in real tridiagonal symmetric matrices is $2$.

Finally, we show how to extend the proof to cover any positive bandwidth. For general real symmetric banded matrices, we still need $n$ parameters for the eigenvalues. To understand the freedom in picking the eigenvectors, observe that (except for degenerate cases happening in Zariski closed proper subsets, which we can safely ignore) the left eigendecomposition $V A = D V$, together with the fact that $A$ is banded, yields the block Krylov relation
\begin{equation}  \label{eq:kryl2}
\mbox{span}(V(:,1:kj+p)) = \mbox{span}(v_1,\dots,v_k,Dv_1,\ldots, D v_k,
\ldots, D^{j-1}v_k, D^j v_1, \ldots, D^j v_p),   
\end{equation}
where $0 \leq p <k$, $kj+p \leq n$, and $v_1, \dots v_k$ are the first $k$ column of $V$, or equivalently the matrix of first $k$ elements of the right eigenvectors of $A$. By equation \eqref{eq:kryl2}, once the orthonormal matrix whose rows are $v_1,\dots,v_k$ has been fixed, everything else is determined (up to signs and except in degenerate cases). Hence, taking into account normalization and orthogonality constraints, one obtains, as required, the counting

$$ n + \sum_{j=1}^{k} (n-j) = \frac{1}{2}(k+1)(2n - k).$$

When a multiple eigenvalue occurs, and still ignoring subtleties that do not matter for the dimensional analysis, the counting modifies as follows. There are $n-1$ degrees of freedom for the eigenvalues. We then get to pick the first $k$ components of the eigenvectors, i.e., $nk$ entries, from which we need to subtract $k$ normalization constraints, $k(k-1)/2$ orthogonality constraints, and $1$ further degree of freedom coming from the redundancy of multiplying the first $k$ rows of the $n \times 2$ matrix representing the eigenspace of the double eigenvalue by any $2 \times 2$ orthogonal matrix. Therefore the dimension of the manifold of symmetric matrices with bandwidth $k$ and multiple eigenvalues is

$$ n-1+nk-k(k+1)/2-1 = \frac{(k+1)(2n-k)}{2} - 2 $$
and we conclude that the codimension is $2$.
\end{proof}

We note that the case $k=0$ of diagonal matrices is exceptional, in that the codimension is $1$ since there is no freedom at all in the eigenvectors, so the only ``loss'' of dimensionality is by forcing two eigenvalues to coincide. Hence,

\begin{proposition}\label{thm:diagonal}
The codimension of derogatory real diagonal matrices, embedded in the manifold of real diagonal matrices, is $1$.
\end{proposition}

\subsection{Hermitian banded matrices}

The case of Hermitian banded matrices is similar to its real counterpart. For bandwidth $0$, there is no distinction between Hermitian and real diagonal matrices. Thus, Proposition~\ref{thm:diagonal} can be applied.

For bandwidths $\geq$ 1, analogously to the real symmetric case, we obtain that the codimension is the same as in the case of bandwidth $n-1$, i.e., full matrices.

\begin{theorem}\label{thm:bandherm}
The real codimension of derogatory banded Hermitian matrices of bandwidth $k > 0$, embedded in the real manifold of Hermitian symmetric matrices of bandwidth $k$, is $3$.
\end{theorem}

\begin{proof}
We sketch the proof, which is similar to the banded real symmetric case. Let $A \in \C^{n \times n}$ be complex Hermitian of bandwidth $k$. The total dimension of this manifold (in real parameters) is clearly $n$ for the diagonals, $2(n-1)$ for the first superdiagonals, then $2(n-2)$, and so on for a total of
$$n + 2 \sum_{j=1}^k (n-j) = n(2k+1) - k(k+1).$$
An equivalent counting is based on the eigendecomposition. There are $n$ real degrees of freedom associated with the eigenvalues. As for the real symmetric case, we can count the degrees of freedom for the eigenvectors by looking at left eigenvectors and at the Krylov process \eqref{eq:kryl2}. (Note that in the Hermitian case the $j$th left eigenvector is the \emph{complex conjugate} of the vector of the $j$th components of right eigenvectors.) Prescribing the eigenspaces is tantamount to choosing a $k \times n$ submatrix of a unitary matrix, quotiented by the choice of $k$ redundant phases. 
\begin{remark}
Discounting phases becomes very subtle with this ``left eigenvectors'' approach, therefore we give more details here. The redundancy is in the independent phases of the right eigenvectors, or equivalently in the independent phases of each given component of the $k$ prescribed left eigenvectors -- this subtracts $n$ real degrees of freedom. However, unlike in the real case, the matrix $A$ is not uniquely determined by the Krylov process \eqref{eq:kryl2}, as there are $n-k$ choices of the phases for the unitary basis of the Krylov subspaces: note that unlike the real case the phases are allowed to vary in a one-real-dimensional interval, and are not just zero-real-dimensional binary choices. Therefore, overall we discard $n-(n-k)=k$ real parameters due to phase redundancy.
\end{remark}
 In terms of real parameters, the argument above gives $k(2n-k) - k = 2kn - k(k+1)$ real degrees of freedom for the eigenvectors. Thus, we recover that the total real dimension is $n(2k+1)-k(k+1)$.

Now we count the dimension of Hermitian banded matrices with multiple eigenvalues. Similarly to the real symmetric case, there is a loss of $1$ degree of freedom for the eigenvalues which can now be described with $n-1$ real parameters. For the eigenvectors, after having picked a $k \times n$ submatrix of a unitary matrix ($k(2n-k)$ real parameters), we now subtract $k+2$ degrees of freedom for the redundancy of the representation. Indeed, this time the counting is as follows: there are $n+2$ redundant parameters in picking a representation for the first component of the right eigenvectors ($n-2$ for simple eigenspaces and $4$ for the double eigenspace). However, the Krylov process \eqref{eq:kryl2} does not uniquely determine $A$ as there are $n-k$ choices of phases; in total we subtract $k+2$ redundancy constraints. In conclusion the dimension is
$$ n-1 + k(2n-k)-(k-2)=n(2k+1)-k(k+1)-3.$$
In other words, the codimension is $3$.
\end{proof}

\begin{figure}[htpb]
  \begin{minipage}[t]{0.5\hsize}
      \includegraphics[width=.9\textwidth]{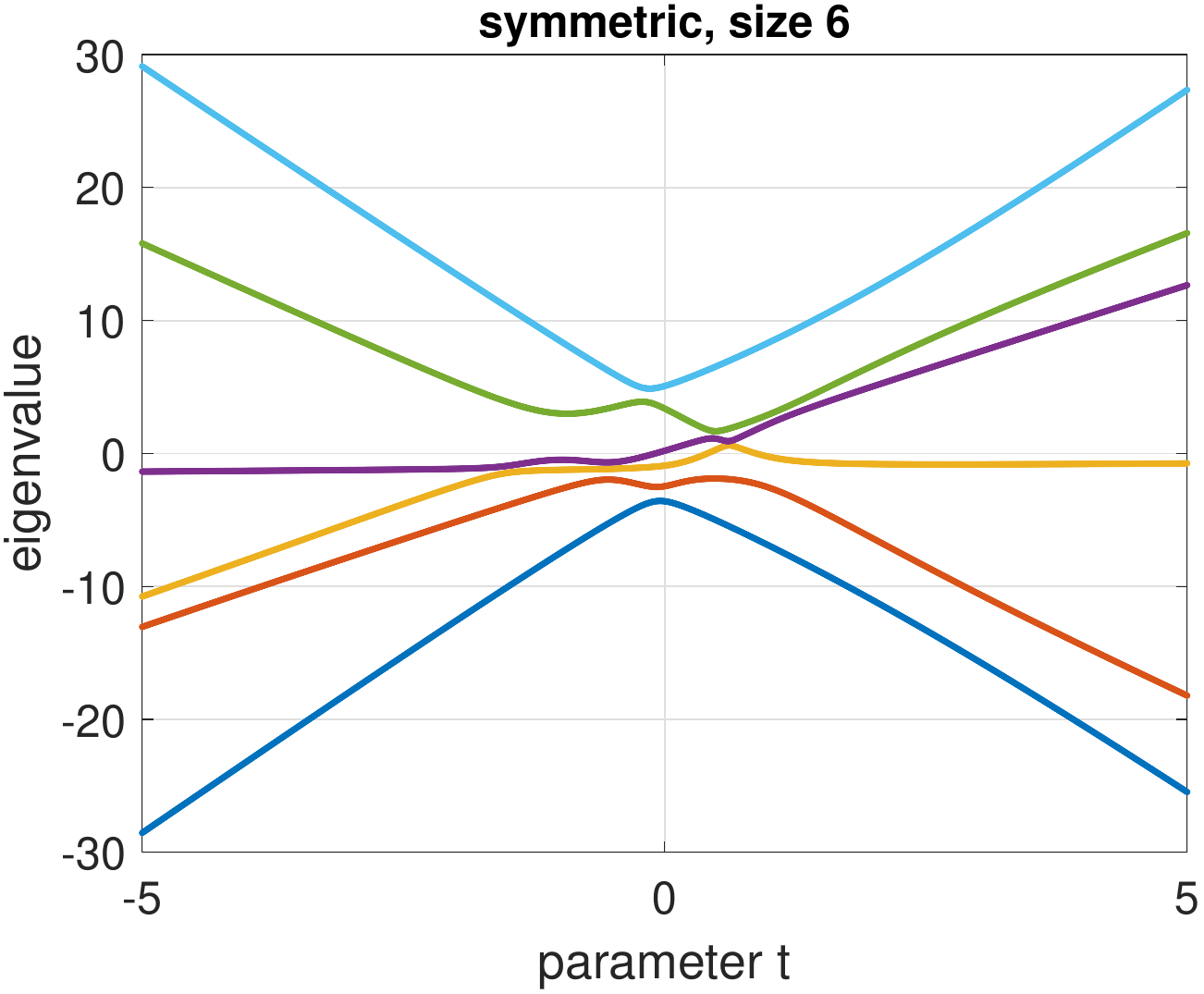}
  \end{minipage}   
  \begin{minipage}[t]{0.5\hsize}
      \includegraphics[width=.9\textwidth]{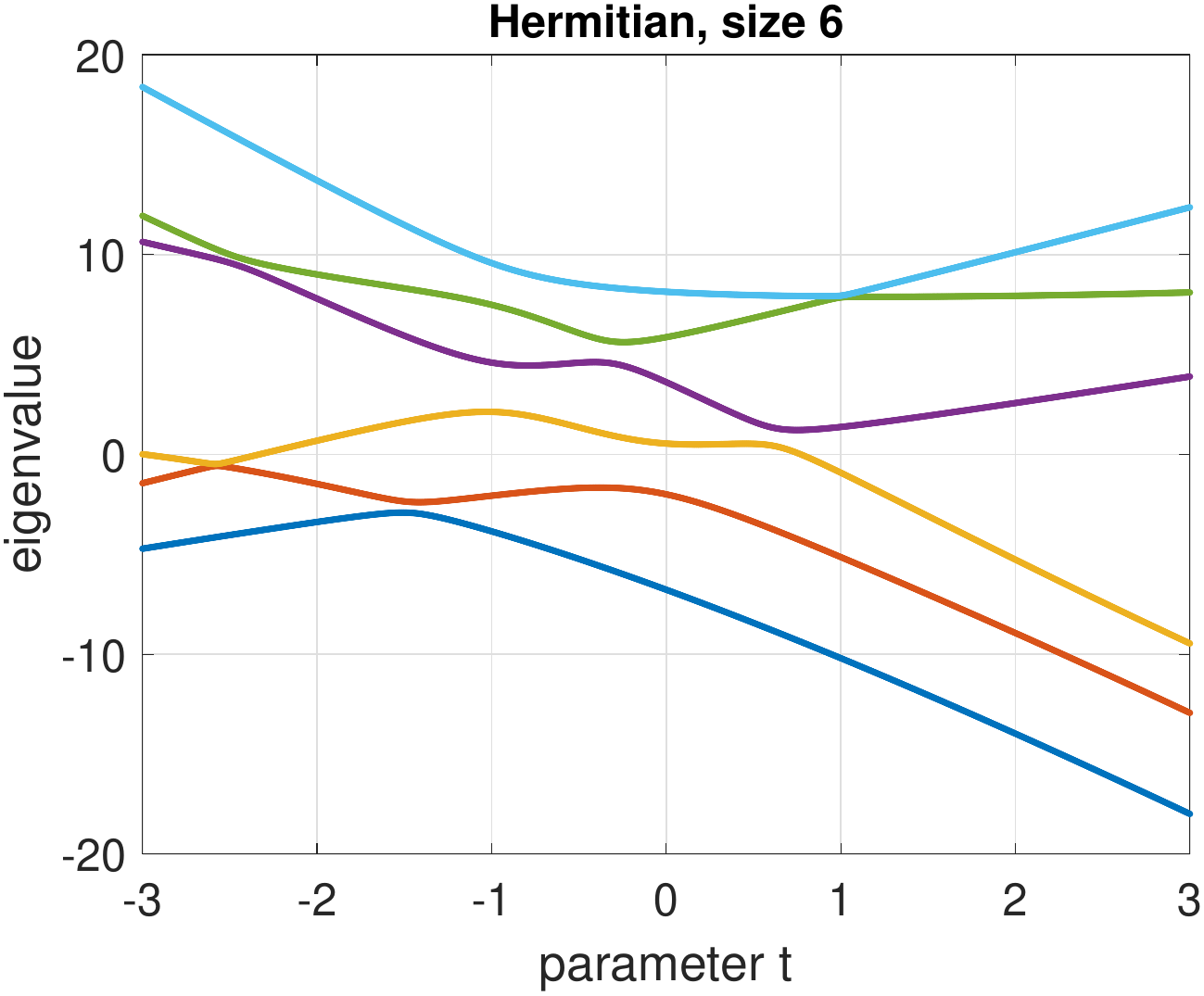}
  \end{minipage}
  \caption{Tridiagonal matrices $A+tB$, codimension 2. Symmetric (left) and Hermitian (right), Codimension 3.}
  \label{fig:trisym}
\end{figure}

\subsection{Banded skew-symmetric matrices}
\begin{theorem}
The real codimension of derogatory banded $n\times n$ skew-symmetric matrices of bandwidth $k > 0$, embedded in the real manifold of skew-symmetric matrices of bandwidth $k$, is $3$when $n$ is odd and it is $1$ when $n$ is even.
\end{theorem}
\begin{proof}
First observe that much of the discussions leading to \eqref{eq:kryl} and \eqref{eq:kryl2} carry over to the skew-symmetric case.

As before, we split the discussion based on the parity of $n$. 
We first discuss even $n$. 

In the real skew-symmetric tridiagonal case, we have only $n-1$ degrees of freedom in the matrix. 
They 
 are accounted for by the eigenvalues, which come in complex conjugate pairs: $n/2$, and eigenvector, $n/2-1$. 
To see this, we use the real block diagonal form of skew-symmetric matrices 
  \begin{equation}
    \label{eq:decompskewsym}
\widetilde VA = \widetilde D\widetilde V , \qquad 
\widetilde D = \mbox{diag}(\big[\begin{smallmatrix}0& a_1\\  -a_1& 0\end{smallmatrix}\big],\big[\begin{smallmatrix}0& a_2\\  -a_2& 0\end{smallmatrix}\big],\ldots,\big[\begin{smallmatrix}0& a_{\frac{n}{2}}\\  -a_{\frac{n}{2}}& 0\end{smallmatrix}\big]),
  \end{equation}
where $a_i\in\mathbb{R}$ 
and $\widetilde V$ is real orthogonal. Now arguing as before, the freedom in choosing $\widetilde V$ reduces to that of choosing the first column $\widetilde v_1$. Note that for each $2\times 2$ block $\widetilde D_i = \big[\begin{smallmatrix}0& a_i\\  -a_i& 0\end{smallmatrix}\big]$ and any orthogonal $2\times 2$ matrix $Q$,  $Q^T\widetilde D_iQ$ is equal to either $\widetilde D_i$ or $-\widetilde D_i$, both of which are allowed in the $\widetilde VA = \widetilde D\widetilde V$ decomposition. Thus for each pair of elements in $\widetilde v_1$ there is one degree of freedom,
hence with the additional constraint $\|\widetilde v_1\|=1$, we see that $\widetilde v_1$ has $n/2-1$ degrees of freedom.

In the derogatory case 
the degrees of freedom are $n/2-1$ in the eigenvalues, and in the eigenvector $n/2-2$.   
Therefore the real codimension is 2.

Now for the general bandwidth, we again start from the real decomposition
$\widetilde VA = \widetilde D\widetilde V$. The freedom in choosing $\widetilde V$ lies in choosing the first $k$ columns of the orthogonal matrix $\widetilde V$, minus the redundancy due to the structure of $\widetilde D$, which is $n/2$. 
Therefore the total degrees of freedom is 
\begin{equation}  \label{eq:dofskewband}
\frac{n}{2}+ \sum_{j=1}^k(n-j) -\frac{n}{2} =\sum_{j=1}^k(n-j)= \frac{k(2n-1-k)}{2},   
\end{equation}
as expected. Now when the matrix is derogatory, as in Section~\ref{sec:evenskew} we split the discussion depending on whether the multiple eigenvalue is zero. When it is zero, the only constraint is that $a_i=0$ for some $i$, and the codimension is 1. When a nonzero eigenvalue $\pm ia$ is multiple, we have $n/2-1$ degrees of freedom in the eigenvalues, and for the eigenvectors, relative to $\sum_{j=1}^k(n-j) -\frac{n}{2}$ in~\eqref{eq:dofskewband}, we further subtract two degrees of freedom, corresponding to the (blockwise) rotation in the equal diagonal blocks $\widetilde D_i=\widetilde D_j$, of the form $Q_i\otimes Q_j$ where $Q_i,Q_j$ are $2\times 2$ orthogonal. It follows that the codimension is 3.

When $n$ is odd, the decomposition~\eqref{eq:decompskewsym} holds with $(n-1)/2$ blocks of size $2\times 2$, and one zero block of size $1\times 1$. The degrees of freedom for bandwidth $k$ is $(n-1)/2$ eigenvalues plus 
the freedom in choosing the first $k$ columns of $\widetilde V$, which is 
$\sum_{j=1}^k(n-j)$, minus the redundancy due to the structure of $\widetilde D$, which is $(n-1)/2$. 
Therefore the total degrees of freedom is 
$\frac{(n-1)}{2}+ \sum_{j=1}^k(n-j) -\frac{(n-1)}{2} =\sum_{j=1}^k(n-j)= \frac{k(2n-1-k)}{2}$, analogous to~\eqref{eq:dofskewband}. 
For derogatory matrices with a nonzero multiple eigenvalue, we subtract one degree for the eigenvalues. For the eigenvectors, the argument is similar to the case where $n$ is even. When zero is a multiple eigenvalue (of multiplicity three), the redundancy comes from the $3\times 3$ zero block in~\eqref{eq:dofskewband}, which involves three degrees of freedom (minus one, which was already subtracted before). The overall codimension is therefore 3. 
\end{proof}

\subsection{Skew-Hermitian banded matrices}
A matrix $A$ is skew-Hermitian and has bandwidth $k$ if and only if $iA$ is Hermitian and has bandwidth $k$. Hence, from Theorem~\ref{thm:bandherm} we immediately have the following result.

\begin{theorem}
The real codimension of derogatory banded skew-Hermitian matrices of bandwidth $k > 0$, embedded in the real manifold of skew-Hermitian symmetric matrices of bandwidth $k$, is $3$.
\end{theorem}

\begin{figure}[htpb]
  \begin{minipage}[t]{0.5\hsize}
      \includegraphics[width=.95\textwidth]{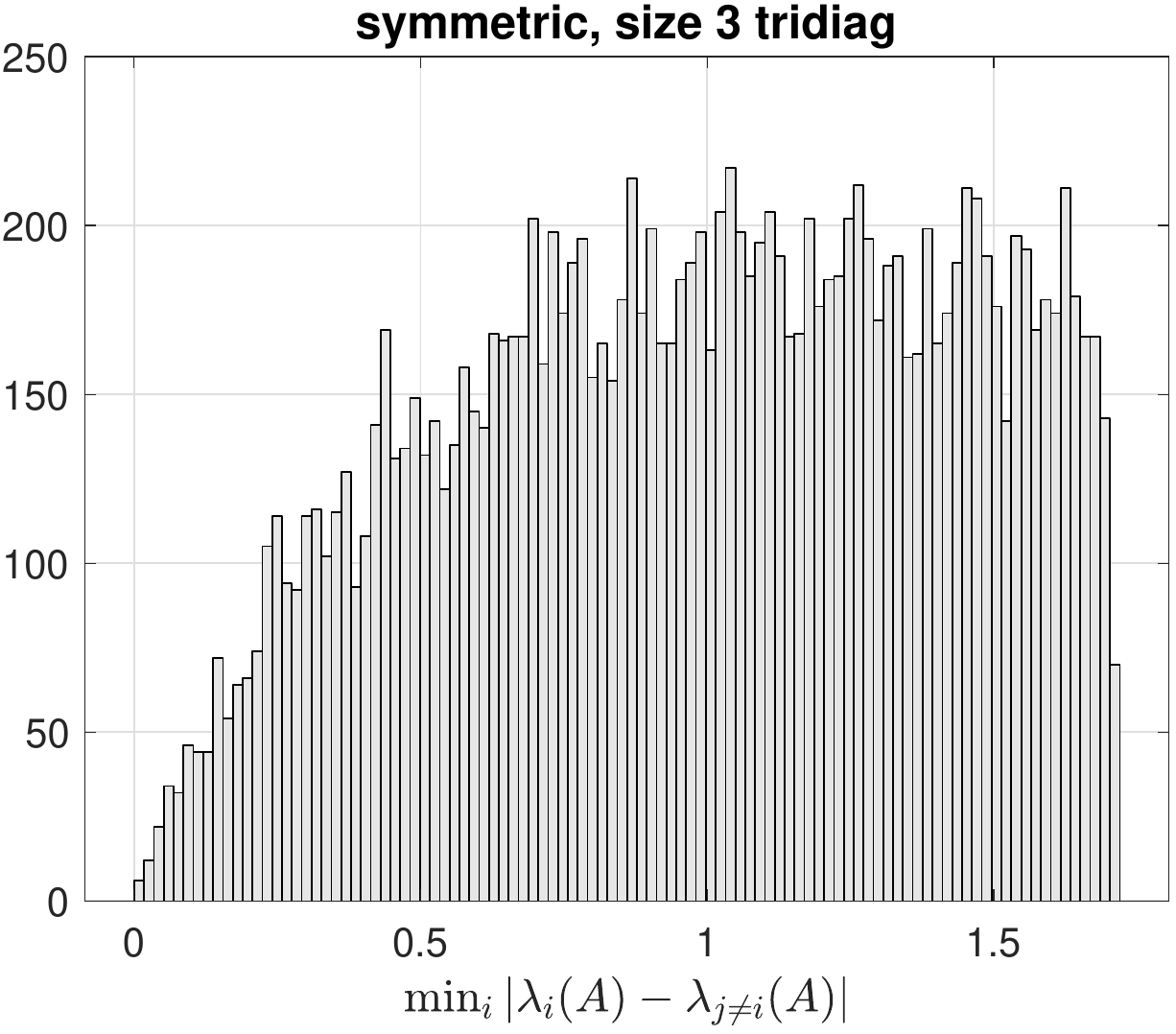}
  \end{minipage}   
  \begin{minipage}[t]{0.5\hsize}
      \includegraphics[width=.95\textwidth]{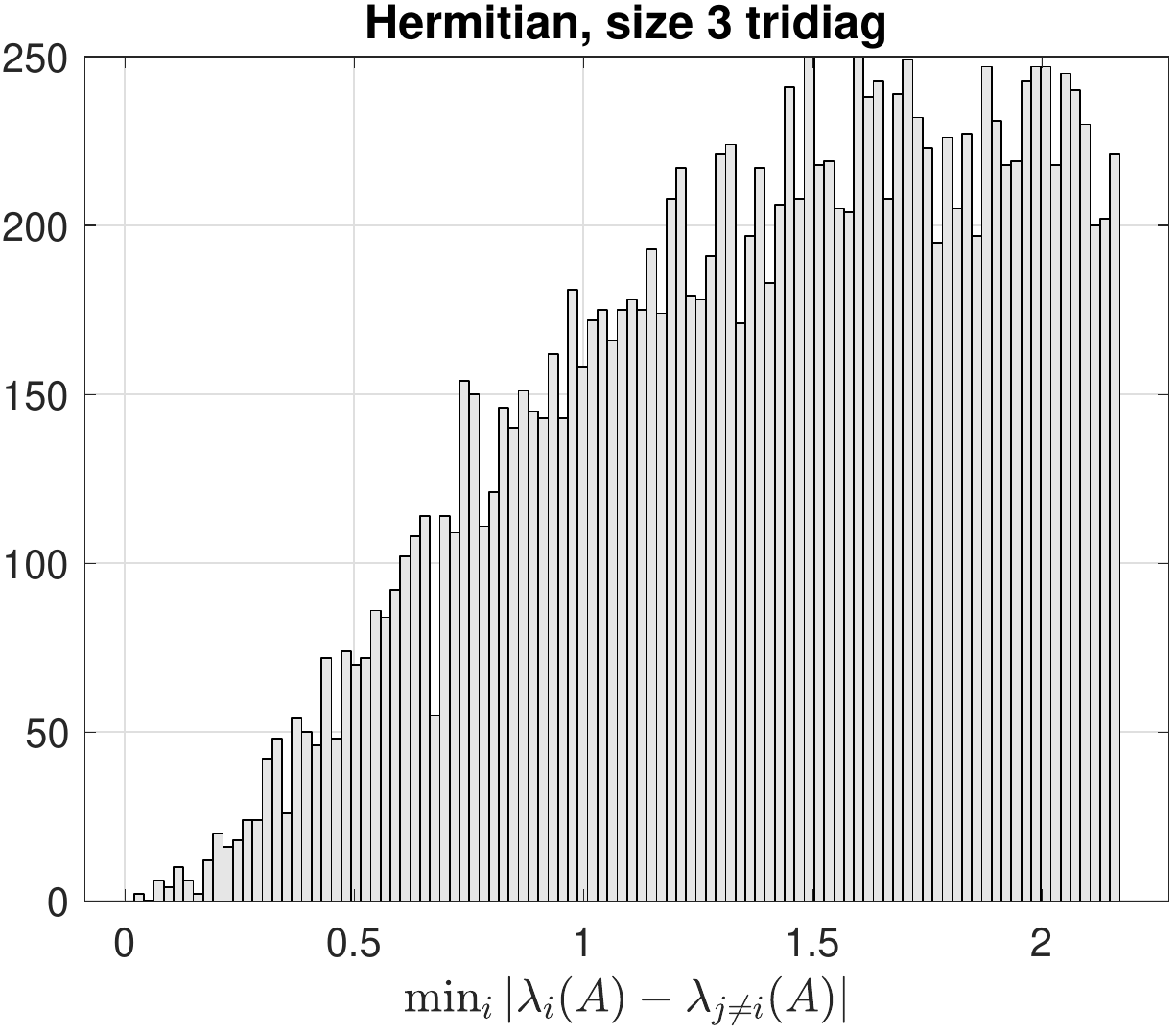}
  \end{minipage}
  \caption{Histogram for $\min_i|\lambda_i(A)-\lambda_{j\neq i}(A)|$, $3\times 3$ tridiagonal matrices (sampled $10^4$ times). Symmetric (left, codimension 2) and Hermitian (right, codimension 3).}
  \label{fig:histtrisym3}
\end{figure}

\begin{figure}[htpb]
  \begin{minipage}[t]{0.5\hsize}
      \includegraphics[width=.95\textwidth]{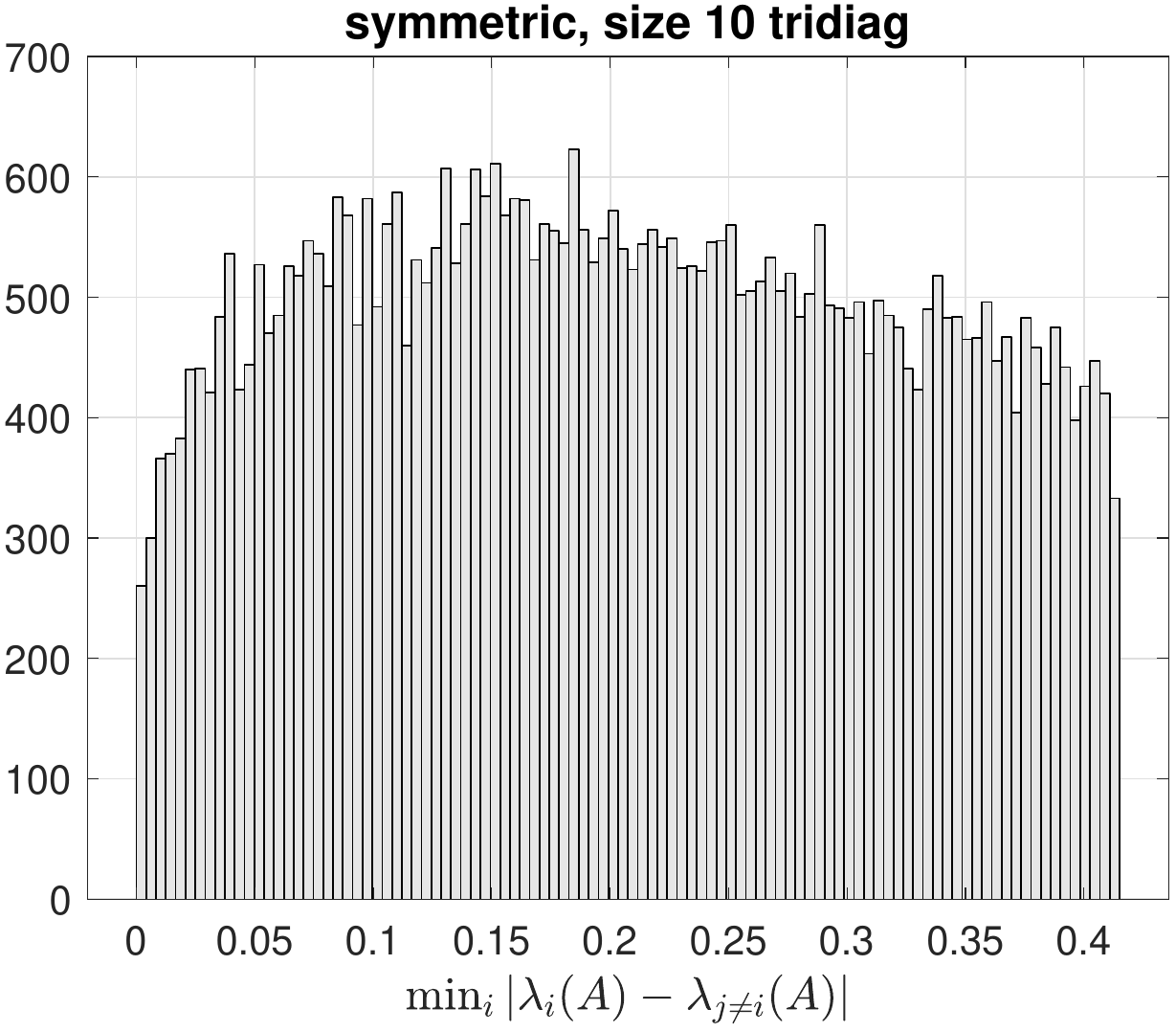}
  \end{minipage}   
  \begin{minipage}[t]{0.5\hsize}
      \includegraphics[width=.95\textwidth]{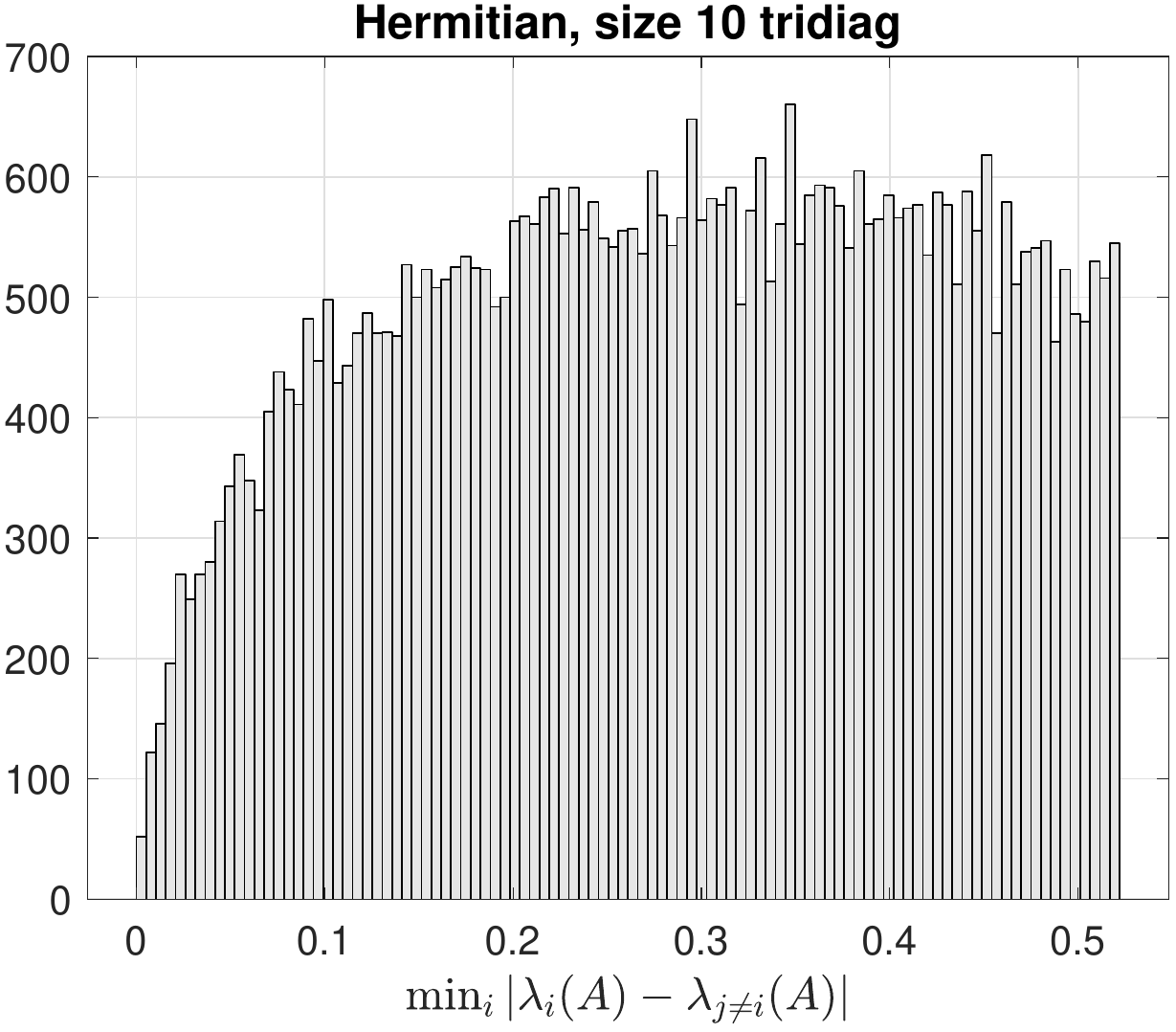}
 \end{minipage}
  \caption{Histogram for $\min_i|\lambda_i(A)-\lambda_{j\neq i}(A)|$, $10\times 10$ tridiagonal matrices (sampled $10^4$ times). Symmetric (left, codimension 2) and Hermitian (right, codimension 3).}
  \label{fig:histtrisym10}
\end{figure}

Experiments suggest that these results on banded matrices are not necessarily reflected in a histogram plot of the eigenvalue distance, especially when the bandwidth is small: When the matrix size $n$ is small, experiments indicate that the dimension counting is accurately reflected, as illustrated in Figure~\ref{fig:histtrisym3}. However, for large $n$, the histograms look significantly different, see Figure~\ref{fig:histtrisym10}. When the bandwidth is large (say $k>2$), the histograms look similar to those of dense matrices.

\begin{figure}[htpb]
  \begin{minipage}[t]{0.5\hsize}
      \includegraphics[width=.95\textwidth]{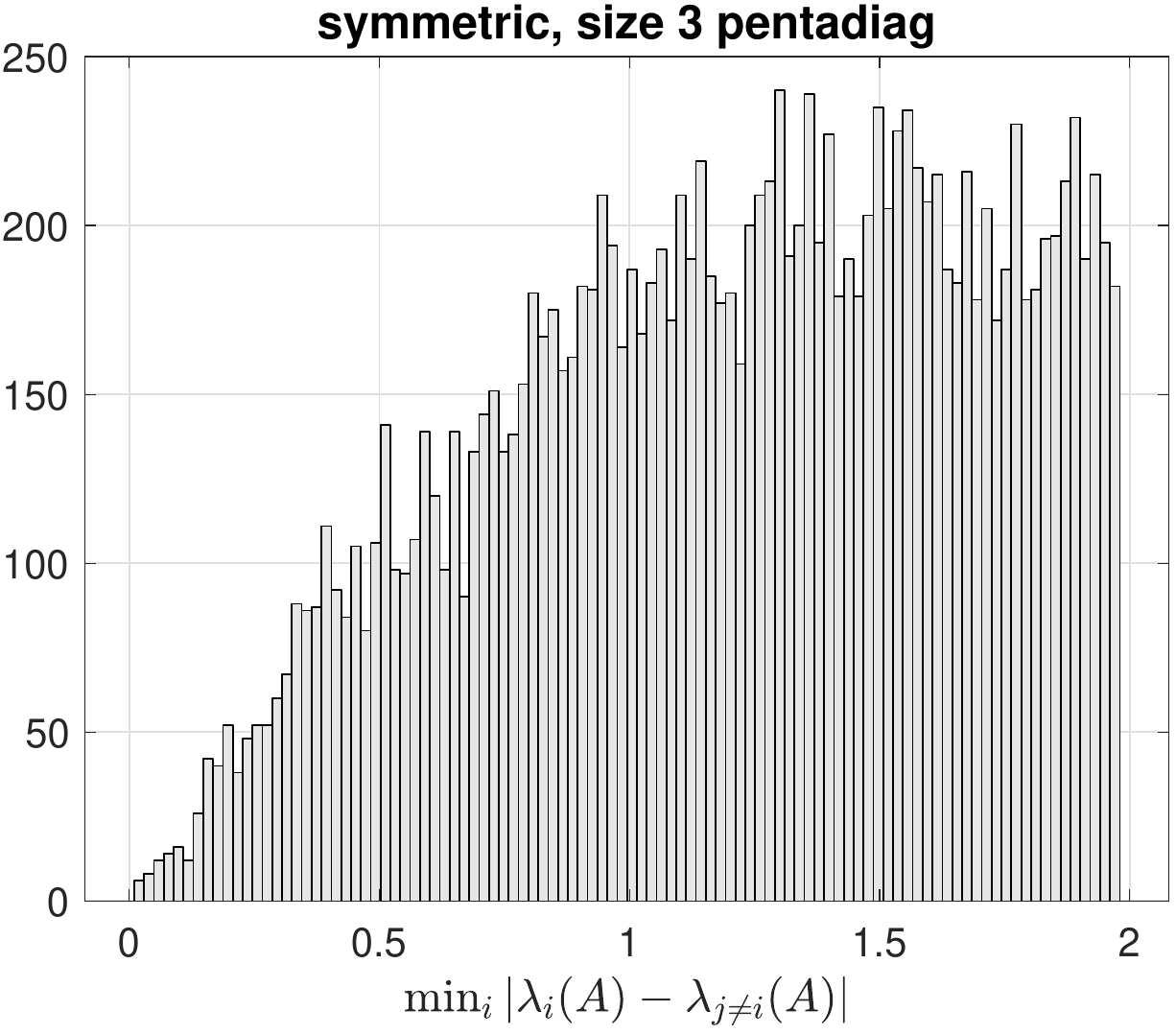}

  \end{minipage}   
  \begin{minipage}[t]{0.5\hsize}
      \includegraphics[width=.95\textwidth]{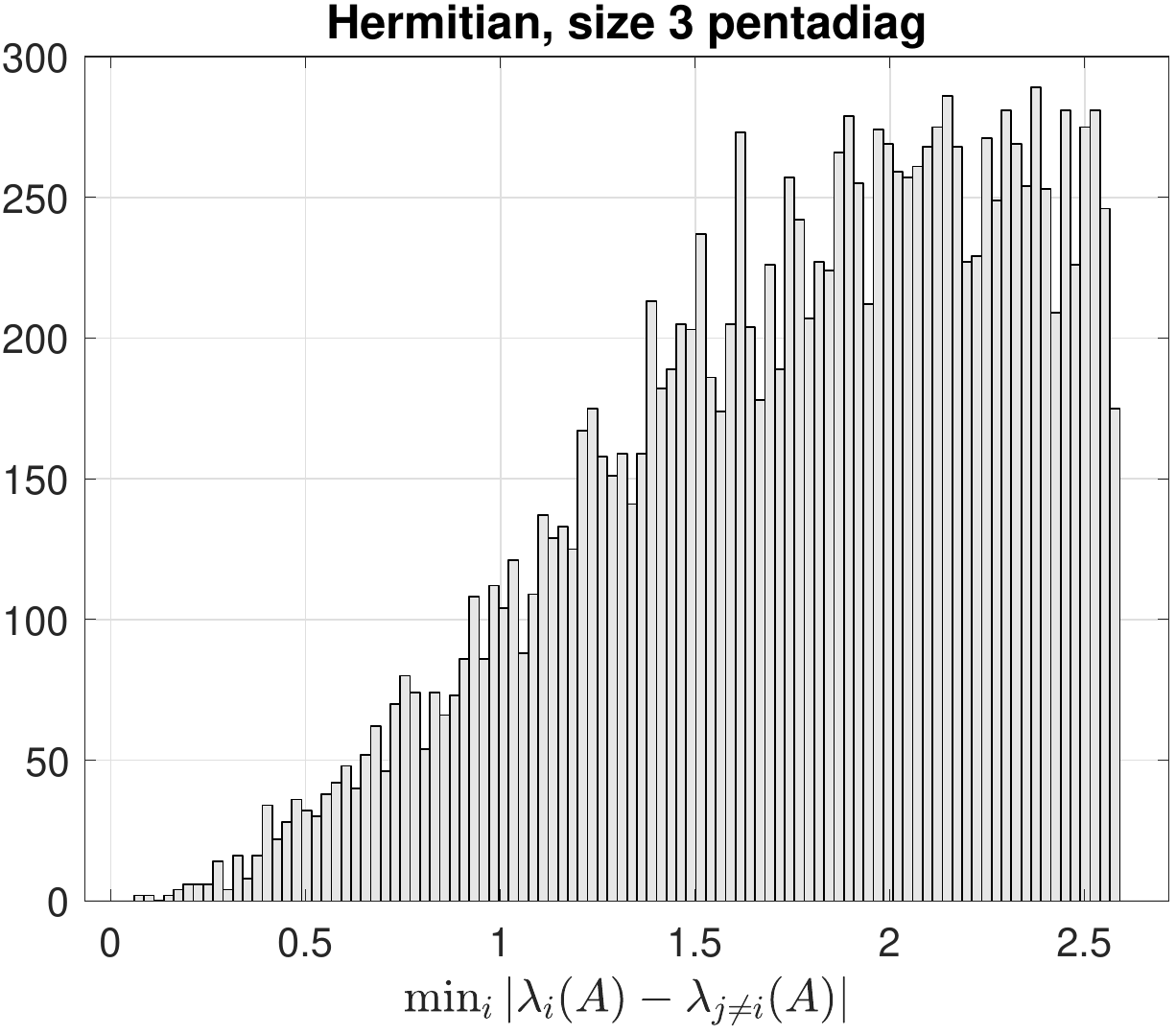}
  \end{minipage}
  \caption{Histogram for $\min_i|\lambda_i(A)-\lambda_{j\neq i}(A)|$, $3\times 3$  pentadiagonal matrices (i.e., dense matrices sampled $10^4$ times). Symmetric (left, codimension 2) and Hermitian (right, codimension 3).}
  \label{fig:histtrisym3ver2}
\end{figure}

\begin{figure}[htpb]
  \begin{minipage}[t]{0.5\hsize}
      \includegraphics[width=.95\textwidth]{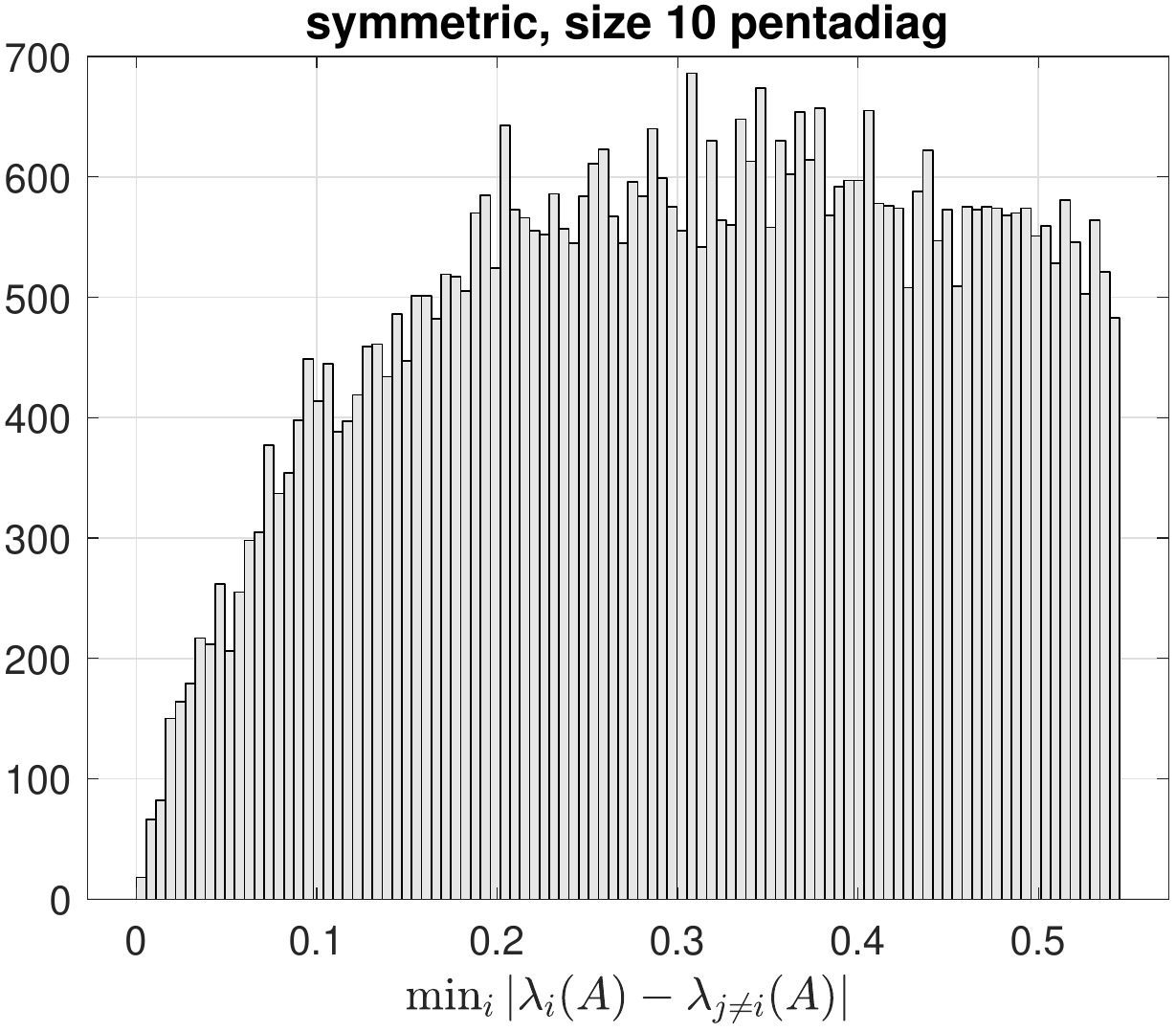}
  \end{minipage}   
  \begin{minipage}[t]{0.5\hsize}
      \includegraphics[width=.95\textwidth]{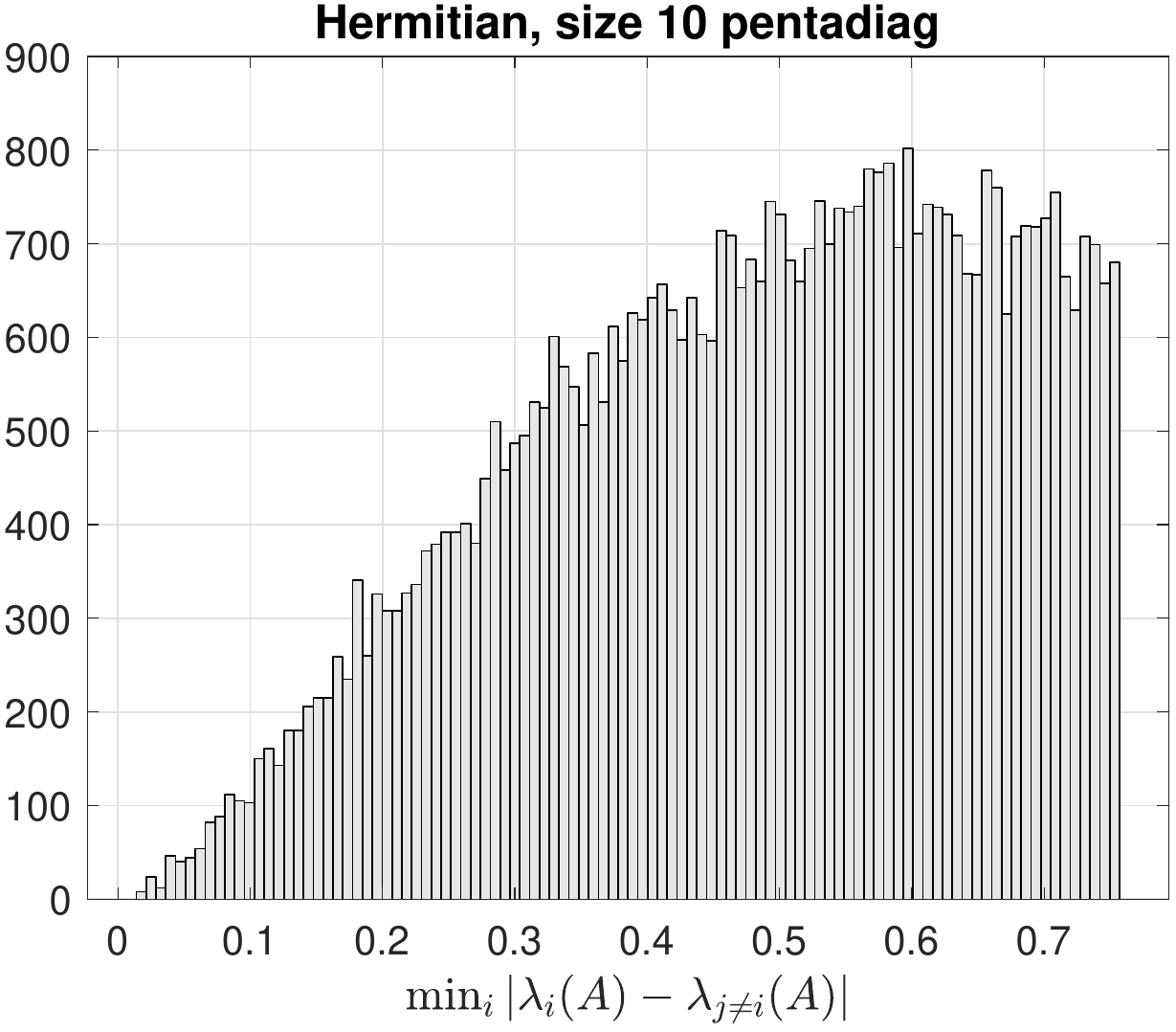}
 \end{minipage}
  \caption{Histogram for $\min_i|\lambda_i(A)-\lambda_{j\neq i}(A)|$, $10\times 10$ pentadiagonal matrices (sampled $10^4$ times). Symmetric (left codimension 2) and Hermitian (right, codimension 3).}
  \label{fig:histtrisym10ver2}
\end{figure}

\section{Singular values}

A singular value decomposition of a matrix $A \in \C^{m\times n}$ is $A = U \Sigma V^*$ where $U \in C^{m \times m}$ and $V \in \C^{n \times n}$ are unitary while $\Sigma \in \R^{m \times n}$ is real, diagonal, and with nonincreasing 
nonnegative diagonal elements \citep{hornjohn}. If $A$ is real, then $U$ and $V$ can be taken real orthogonal. There is a link~(e.g. \cite[\S~I.4.1]{stewart-sun:1990}) between singular values of an unstructured matrix and eigenvalues of symmetric related to $A$, such as $A^*A$ or
$$\begin{bmatrix}
0 & A\\
A^* & 0
\end{bmatrix}.$$
These matrices are however highly structured because of their zero pattern or their positive semidefiniteness. 
In this Section, we study the phenomenon of singular value avoidance.

Keller~\cite[\S~8]{keller2008} showed that the codimension of real $m \times n$ matrices with repeated singular values, embedded in $\R^{m\times n}$, is $2$. This can be obtained using our line of arguments, via the degrees of freedom in the singular vectors and values. We next extend the result to complex matrices, which---as in the symmetric vs. Hermitian contrast---has a larger codimension. 

\begin{theorem}\label{thm:svd}
The real codimension of complex $m \times n$ matrices with repeated singular values, embededd in $\C^{m\times n}$, is $3$.
\end{theorem}

\begin{proof}
We assume without loss of generality that $m \geq n$. As usual for the complex case, we will need particular care for the redundancy of phases. We first recount the real dimension of $m \times n$ complex matrices via their singular value decomposition. There are $n$ degrees of freedom for the singular values, $n^2$ for the choice of an $n\times n$ unitary matrix, $2mn-n^2$ degrees of freedom to choose the first $n$ columns of an $m \times m$ unitary matrix, and we then must subtract $n$ redundant parameters for the phases (\emph{not} $2n$ because the ratio of the phases of the right and left singular vectors \emph{is} uniquely determined by $A$). We thus obtain the expected value of $2mn$ for the real dimension.

For matrices with repeated singular value, we need $n-1$ real parameters for the singular values, $n^2-4$ for the nondegenerate right singular spaces, 
$(2m-n+2)(n-2)$
for the nondegenerate left singular spaces, $4$ for the multiple right singular space, 
$4(m-n+1)$
 for the multiple left singular space, minus $n-2$ redundant phases for the nondegenerate singular spaces, minus $4$ real parameters for the degenerate singular spaces (we can multiply the columns reprenseting a basis for the left and right singular spaces times the same $2 \times 2$ unitary matrix). In total, the dimension is
$$
(n-1)+(n^2-4)+ (2m-n+2)(n-2)+4+4(m-n+1)-(n-2)-4 = 2mn-3.
$$
\end{proof}

\noindent Finally, we note that structured matrices can easily have multiple singular values. For example, a symmetric matrix has a multiple singular value when it has an eigenvalue pair $\pm \lambda$; the codimension is 1. As a more extreme example, the singular values of orthogonal matrices are clearly always all 1. 

\section{Conclusions}
We have studied the codimension of submanifolds of derogatory structured matrices, thus predicting whether the phenomenon of eigenvalue avoidance is to be expected or not. Interesting potential continuations of this line of research include the following:
(1) Extending from Hermitian to other structures also the investigation of near-crossing, i.e., the fact that eigenvalues, even if they do not intersect, tend to come very close to each other~\cite{betcke2004computations}; 
(2) Investigating more structure and possibly extending this work to structured matrix pencils and matrix polynomials.

\section*{Acknowledgements}

We thank Nick Trefethen for introducing us to this subject during a lunch at Balliol College in Oxford.

\bibliographystyle{abbrv}
\bibliography{bib}

\begin{thebibliography}{10}

\bibitem{arnold2012geometrical}
V.~I. Arnol'd.
\newblock {\em Geometrical Methods in the Theory of Ordinary Differential
  Equations}, volume 250.
\newblock Springer Science \& Business Media, 2012.

\bibitem{betcke2004computations}
T.~Betcke and L.~N. Trefethen.
\newblock Computations of eigenvalue avoidance in planar domains.
\newblock {\em PAMM}, 4(1):634--635, 2004.

\bibitem{bhat:96}
R.~Bhatia.
\newblock {\em Matrix Analysis}.
\newblock Graduate Texts in Mathematics, vol. 169. Springer, New York, 1996.

\bibitem{dieci1999smooth}
L.~Dieci and T.~Eirola.
\newblock On smooth decompositions of matrices.
\newblock {\em SIAM J. Matrix Anal. Appl.}, 20(3):800--819, 1999.

\bibitem{DPP18}
L.~Dieci, A.~Papini, and A.~Pugliese.
\newblock Coalescing points for eigenvalues of banded matrices depending on
  parameters with application to banded random matrix functions.
\newblock {\em Numer. Algorithms}, 80:1241--1266, 2019.

\bibitem{edel98}
A.~Edelman, T.~A. Arias, and S.~T. Smith.
\newblock {The geometry of algorithms with orthogonality constraints}.
\newblock {\em SIAM J. Matrix Anal. Appl.}, {20}({2}):{303--353}, {1998}.

\bibitem{golubbook4th}
G.~H. Golub and C.~F. Van~Loan.
\newblock {\em Matrix {C}omputations}.
\newblock {The Johns Hopkins University Press}, 4th edition, 2012.

\bibitem{hornjohn}
R.~A. Horn and C.~R. Johnson.
\newblock {\em Matrix {A}nalysis}.
\newblock Cambridge University Press, second edition, 2012.

\bibitem{kato}
T.~Kato.
\newblock {\em Perturbation Theory for Linear Operators}.
\newblock Springer-Verlag, 2nd edition, 1966.

\bibitem{keller2008}
J.~B. Keller.
\newblock Multiple eigenvalues.
\newblock {\em Linear Algebra Appl.}, 429(8-9):2209--2220, 2008.

\bibitem{Laxlaa}
P.~Lax.
\newblock {\em Linear Algebra and Its Applications}.
\newblock Wiley, 2nd edition, 2007.

\bibitem{lee}
J.~M. Lee.
\newblock {\em Introduction to Smooth Manifolds}.
\newblock Springer, 2003.

\bibitem{MNTX16}
V.~Mehrmann, V.~Noferini, F.~Tisseur, and H.~Xu.
\newblock On the sign characteristics of {H}ermitian matrix polynomials.
\newblock {\em Linear Algebra Appl.}, 511:328--364, 2016.

\bibitem{rellich37}
F.~Rellich.
\newblock St\"{o}rungstheorie der {S}pektralzerlegung {I}.
\newblock {\em Math. Anal.}, 113:600--619, 1937.

\bibitem{rellichbook}
F.~Rellich.
\newblock {\em Perturbation Theory of Eigenvalue Problems}.
\newblock Gordon and {B}reach, 1969.

\bibitem{smale1958}
S.~Smale.
\newblock Regular curves on riemannian manifolds.
\newblock {\em Trans. Am. Math. Soc.}, 87(2):492--512, 1958.

\bibitem{stewart-sun:1990}
G.~W. Stewart and J.-G. Sun.
\newblock {\em Matrix Perturbation Theory (Computer Science and Scientific
  Computing)}.
\newblock Academic Press, 1990.

\bibitem{teytel1999rare}
M.~Teytel.
\newblock How rare are multiple eigenvalues?
\newblock {\em Comm. Pure Appl. Math.}, 52(8):917--934, 1999.

\bibitem{Thom2}
R.~Thom.
\newblock Quelques proprietes globales des varietes differentiables.
\newblock {\em Comment. Math. Helv.}, 28(1):17--86, 1954.

\bibitem{Thom}
R.~Thom.
\newblock Un lemme sur les applications differentiables.
\newblock {\em Bol. Soc. Mat. Mexicana}, 2(1):59--71, 1956.

\bibitem{uhlenbeck1976generic}
K.~Uhlenbeck.
\newblock Generic properties of eigenfunctions.
\newblock {\em Am. J. Math.}, 98(4):1059--1078, 1976.

\bibitem{von1993verhalten}
J.~von Neumann and E.~P. Wigner.
\newblock {\"U}ber das verhalten von eigenwerten bei adiabatischen prozessen.
\newblock {\em Phys. Z.}, 30:467--470, 1929.

\end{thebibliography}

\end{document}